\documentclass[10pt]{article}
\usepackage[a4paper, total={7in, 9in}]{geometry} 

\usepackage{amsthm}
\usepackage{amsfonts}
\usepackage{amssymb}
\usepackage{amsmath}
\usepackage{mathtools}
\usepackage[maxbibnames=99]{biblatex}
\usepackage{xcolor}
\usepackage{esint} 
\usepackage{tikz-cd}
\usepackage{blindtext}
\usepackage{titlesec}
\usepackage{hyperref}
\usepackage{nomencl}

\newtheorem{innercustomgeneric}{\customgenericname}
\providecommand{\customgenericname}{}
\newcommand{\newcustomtheorem}[2]{%
  \newenvironment{#1}[1]
  {%
   \renewcommand\customgenericname{#2}%
   \renewcommand\theinnercustomgeneric{##1}%
   \innercustomgeneric
  }
  {\endinnercustomgeneric}
}

\newcustomtheorem{customthm}{Theorem}
\newcustomtheorem{customlemma}{Lemma}
\newcustomtheorem{customC}{Corollary}

\newtheorem{thm}{Theorem}[section]
\newtheorem{corollary}[thm]{Corollary}
\newtheorem{lemma}[thm]{Lemma}
\newtheorem{proposition}[thm]{Proposition}

\newtheorem*{thm*}{Theorem}
\newtheorem*{corollary*}{Corollary}
\newtheorem*{lemma*}{Lemma}
\newtheorem*{proposition*}{Proposition}

\theoremstyle{definition}
\newtheorem{definition}[thm]{Definition}
\newtheorem*{definition*}{Definition}
\newtheorem{remark}[thm]{Remark}

\newtheorem*{remark*}{Remark}

\newcommand{\bb}[1]{\mathbb{#1}}
\newcommand{\ssf}[1]{\mathsf{#1}}
\newcommand{\supp}{\mathsf{supp}}
\newcommand{\Per}{\mathsf{Per}}
\newcommand{\lip}{\mathsf{lip}}
\newcommand{\Lip}{\mathsf{Lip}}
\newcommand{\Tan}{\mathsf{Tan}}

\newcommand{\Ch}{\mathsf{Ch}}
\newcommand{\Epi}{\mathsf{Epi}}
\newcommand{\Div}{\mathsf{div}}
\newcommand{\Ric}{\mathsf{Ric}}

\newcommand{\Hess}{\mathsf{Hess}}
\newcommand{\RCD}{\mathsf{RCD}}
\newcommand{\Graph}{\mathsf{Graph}}

\newcommand{\aH}{\mathsf{H}}
\newcommand{\BV}{\mathsf{BV}}
\newcommand{\W}{\mathsf{W}}
\newcommand{\X}{\mathsf{X}}
\newcommand{\G}{\mathsf{G}}
\newcommand{\M}{\mathsf{M}}

\newcommand{\Osc}{\mathsf{Osc}}
\newcommand{\sd}{\mathsf{d}}
\newcommand{\sL}{\mathsf{L}}
\newcommand{\m}{\mathfrak{m}}
\newcommand{\mres}{\mathbin{\vrule height 1.6ex depth 0pt width
0.13ex\vrule height 0.13ex depth 0pt width 1.3ex}}
\newcommand{\mms}[3]{(\mathsf{#1},\mathsf{#2},\mathfrak{#3})}

\makenomenclature 

\renewcommand{\nompreamble}{The symbols below are listed following the order of appearence in the paper.}

\title{Minimal Surface Equation and Bernstein Property on $\RCD$ spaces}
\author{ \footnote{Mathematical Institute, University of Oxford, Oxford, UK. E-mail address: alessandro.cucinotta@maths.ox.ac.uk} Alessandro Cucinotta}

\begin{document}

\maketitle

\begin{abstract}
We show that if $(\X,\sd,\m)$ is an $\RCD(K,N)$ space and $u \in \W^{1,1}_{loc}(\X)$ is a solution of the minimal surface equation, then $u$ is harmonic on its graph (which has a natural metric measure space structure). If $K=0$ this allows to obtain an Harnack inequality for $u$, which in turn implies the Bernstein property, meaning that any positive solution to the minimal surface equation must be constant.
As an application, we obtain oscillation estimates and a Bernstein Theorem for minimal graphs in products $\M \times \bb{R}$, where $\M$ is a smooth manifold (possibly weighted and with boundary) with non-negative Ricci curvature.
\end{abstract}

\medskip
\noindent
\small{
\textbf{Keywords}: Minimal surface equation, Harnack inequality, Ricci curvature.

\smallskip
\noindent
\textbf{Msc2020}: 53A10, 53C21.
}

\tableofcontents

\section{Introduction}
In \cite{BGM} Bombieri, De Giorgi and Miranda showed that the only entire positive solutions of the minimal surface equation in Euclidean space are the constant functions. If we replace the Euclidean space with a Riemannian manifold the validity of the aforementioned result depends on the geometry of the manifold. For this reason the following definition was introduced in \cite{CMMR}.

\begin{definition*}
    We say that a Riemannian manifold $(\M,g)$ has the Bernstein Property if the only entire positive solutions of the minimal surface equation are the constant functions.
\end{definition*}

For example, positive non constant solutions of the minimal surface equation in the hyperbolic plane were constructed in \cite{NR}, while the Bernstein Property holds for manifolds with non-negative Ricci curvature and a lower sectional curvature bound thanks to \cite{RSS}. Recently this was improved in \cite{CMMR} to manifolds with non-negative Ricci curvature and no sectional curvature constraint (see also \cite{cmmr23} for an even stronger result in the same fashion), while Ding proved in \cite{Ding} that the Bernstein Property holds in manifolds that are doubling and support a Poincaré inequality.
\par
The fact that manifolds with non-negative Ricci curvature have the Bernstein Property and the recent generalization of some properties of minimal surfaces to $\RCD$ spaces (i.e. metric measure spaces with a notion of lower bound on the Ricci curvature) suggest that an analogue of the Bernstein Property might hold in this setting as well. 
We recall that an $\RCD(K,N)$ space is a metric measure space where $K \in \bb{R}$ plays the role of a lower bound on the Ricci curvature, while $N \in [1,+\infty)$ plays the role of an upper bound on the dimension. This class includes measured Gromov-Hausdorff limits of smooth manifolds of fixed dimension with uniform Ricci curvature lower bounds and finite dimensional Alexandrov spaces with sectional curvature bounded from below.
With this in mind we can state the main result of this paper.

\begin{customthm}{1} \label{CT1}
    Let $(\X,\sd,\m)$ be an $\RCD(0,N)$ space and let $u \in \W^{1,1}_{loc}(\X)$ be a solution of the minimal surface equation. If $u$ is positive, then it is constant.
\end{customthm}

The previous theorem, as anticipated, is part of a wider class of recent results that aim at generalizing to the non smooth setting properties of minimal surfaces (\cite{APPS}, \cite{FSM}, \cite{BPSrec}, \cite{Weak}, etc.). 
Moreover, specializing Theorem \ref{CT1} to the smooth category, we obtain that the Bernstein Property holds for certain weighted manifolds with boundary. This is the content of Theorem \ref{CT2}. Given a manifold $(\ssf{M},g)$ we denote by $\ssf{Vol}_g$ its volume measure and by $\ssf{d}_g$ its distance. If $V:\ssf{M} \to \bb{R}$ is a smooth function, we say that the metric measure space $(\ssf{M}^n,\ssf{d}_g,e^{-V} \ssf{Vol}_g)$ is a weighted manifold. Given $\Omega \subset \M$, we say that a function $u \in C^{\infty}(\Omega)$ is a solution of the weighted minimal surface equation on $\Omega \setminus \partial \M$ if
\[
\Div \Big( \frac{e^{-V} \nabla u}{\sqrt{1+|\nabla u|^2}} \Big)=0 \quad \text{on } \Omega \setminus \partial \M.
\]
We say that the boundary of a manifold with boundary is convex if its second fundamental form w.r.t. the inward pointing unit normal is positive.

\begin{customthm}{2} \label{CT2}
Let $(\ssf{M}^n,\ssf{d}_g,e^{-V} \ssf{Vol}_g)$ be a weighted manifold with convex boundary such that there exists $N>n$ satisfying
    \begin{equation} \label{CE1}
    \Ric_\M+\Hess_V-\frac{\nabla V \otimes \nabla V}{N-n} \geq 0 \quad \text{ on } \ssf{M} \setminus \partial \ssf{M}.
    \end{equation}
If $u \in C^{\infty}(\M)$ is a positive solution of the weighted minimal surface equation on $\M \setminus  \partial \M$ whose gradient vanishes on $\partial \M$, then $u$ is constant.
\end{customthm}

The previous result is new (to the best of our knowledge) both in the boundaryless weighted setting and in the framework of unweighted manifolds with boundary. Theorem \ref{CT2}, although being a direct consequence of Theorem \ref{CT1}, can also be obtained independently by more straightforward arguments that rely on the smoothness of the space.
\par
A second consequence of Theorem \ref{CT1} is that the oscillation of minimal graphs in an appropriate class of pointed manifolds grows with a uniform rate as one moves away from the base point in each manifold. This is stated precisely in Theorem \ref{CT3} below. 
A similar result for solutions to elliptic problems in Euclidean spaces was obtained in \cite{moserosc}.
Given a pointed metric space $(\M,\sd,x)$, $r>0$ and $f :B_r(x) \to \bb{R}$, we define 
\[
\Osc_{x,r}(f):=\sup \{|f(y)-f(x)|:y \in B_r(x) \},
\]

\begin{customthm}{3} \label{CT3}
    Let $n \in \bb{N}$ be fixed. For every $T,t,r>0$ with $T>t$ there exists $R>r$ such that if $(\M^n,g,x)$ is a pointed manifold with non-negative Ricci curvature and $u \in C^{\infty}(B_{R}(x))$ is a solution of the minimal surface equation such that $\Osc_{x,r}(u) \geq t$, then
     $
    \Osc_{x,R}(u) \geq T
    $.
\end{customthm}

In Section \ref{S3} we actually prove a more general result involving weighted manifolds with boundary (and even $\RCD(0,N)$ spaces, see Theorem \ref{theorem3stronger} and Corollary \ref{C12}) but we preferred to state Theorem \ref{CT3} in this form for simplicity. 
These results essentially follow by combining the argument of \cite[Remark $4.5$]{Ding} and the Harnack inequality (see Theorem \ref{T10}) that is at the core of Theorem \ref{CT1}. 
\par 
We now turn our attention to the proof of Theorem \ref{CT1}. To this aim let $(\X,\sd,\m)$ be an $\RCD(0,N)$ space and let $u \in \W^{1,1}_{loc}(\X)$ be a solution of the minimal surface equation. 
The proof that we give follows the one in \cite{Ding} to prove the analogous result for manifolds supporting a doubling volume measure and a Poincaré inequality. In that work, Ding shows first that Sobolev and Poincaré inequalities for the height function hold on the graph of $u$. Then, he uses a Moser-type iteration argument, which relies on the fact that $u$ is harmonic on its graph, to obtain the Harnack inequality for $u$ and the Bernstein Property. \par 
The first obstacle for adapting the previously outlined strategy is that the graph of $u$ in our setting has very little structure, and it is not clear what it means for $u$ to be harmonic on its graph. A key result in order to give the graph of $u$ a metric measure space structure is Theorem \ref{CT4}. We denote by $\Epi(u)$ the epigraph of $u$.

\begin{customthm}{4} \label{CT4}
    Let $(\X,\sd,\m)$ be an $\RCD(K,N)$ space, let $\Omega \subset \X$ be open and let $u \in \W^{1,1}_{loc}(\Omega)$. The following conditions are equivalent.
    \begin{enumerate}
        \item $\Epi(u)$ is locally perimeter minimizing in $\X \times \bb{R}$.
        \item $\Epi(u)$ is perimeter minimizing in $\X \times \bb{R}$.
        \item $u$ solves the minimal surface equation on $\X$.
    \end{enumerate}
\end{customthm}

\begin{remark*}
    Thanks to Theorem \ref{CT4}, if we assume that the spaces in question satisfy appropriate parabolicity constraints, then Theorems \ref{CT1}, \ref{CT2} and \ref{CT3} follow from the fact that parabolic $\RCD(0,N)$ spaces have the Half Space Property (see \cite{C1}). We stress that the proofs of these theorems in our setting (i.e. without any parabolicity assumption) use completely different techniques from the ones of \cite{C1}.
\end{remark*}

Thanks to the previous theorem, if $u \in \W^{1,1}_{loc}(\X)$ is a solution of the minimal surface equation, we can consider the closed representative of $\Epi(u)$ (which exists since this is a perimeter minimizing set), and we can define a complete separable metric measure space which plays the role of the graph of $u$. More precisely, we consider $(\G(u),\sd_g,\m_g)$, where $\G(u):=\partial \Epi(u)$, $\sd_g$ is the restriction of the product distance in $\X \times \bb{R}$ to $\G(u)$ and $\m_g$ is the restriction of the perimeter measure of $\Epi(u)$ to $\G(u)$. We also denote by $u_g:\G(u) \to \bb{R}$ the height function on $\G(u)$.
With this definition of graph space, the Sobolev and Poincaré inequalities for $u_g$ on $\G(u)$ follow mimicking the proofs in \cite{Ding} (with the due modifications), so that the problem reduces to proving that $u_g$ is harmonic on $\G(u)$ in a generalized sense that allows to repeat the iteration argument in \cite{Ding}. \par 
The key results in this sense are given by Theorems \ref{CT5} and \ref{CT6}. Given $f:\G(u) \to \bb{R}$ and $x \in \G(u)$, we denote by $\lip_g(f)(x)$ the local Lipschitz constant of $f$ at $x$ (see \eqref{E66}). At every point of $\G(u)$ where it makes sense, given a second function $g:\G(u) \to \bb{R}$, we define
    \[
    \lip_g(f) \cdot \lip_g(g):=\frac{1}{4}(\lip_g(f+g)^2-\lip_g(f-g)^2).
    \]
    We denote by $\Lip(\G(u))$ and $\Lip_c(\G(u))$ respectively the set of Lipschitz continuous and compactly supported Lipschitz continuous functios on $\G(u)$.
Theorem \ref{CT5} shows that the previously defined product, when restricted to Lipschitz functions, behaves according to the usual rules of products of gradients.

\begin{customthm}{5} \label{CT5}
    Let $(\X,\sd,\m)$ be an $\RCD(K,N)$ space, let $\Omega \subset \X$ be open, and let $u \in \W^{1,1}_{loc}(\Omega)$ be a solution of the minimal surface equation. The function $\cdot:\Lip(\G(u)) \times \Lip(\G(u)) \to \sL^\infty(\G(u))$ given by 
    \[
    (\phi_1,\phi_2) \mapsto \lip_g(\phi_1) \cdot \lip_g(\phi_2)
    \]
    is symmetric, bilinear, it satisfies the chain rule and the Leibniz rule in both entries and $\lip_g (\phi_1) \cdot \lip_g(\phi_2) \leq \lip_g(\phi_1) \lip_g(\phi_2)$.
\end{customthm}
Theorem \ref{CT6} is then the analogue in our setting of the fact that in the smooth category $u_g$ would be harmonic on its graph. 

\begin{customthm}{6} \label{CT6}
Let $(\X,\sd,\m)$ be an $\RCD(K,N)$ space, let $\Omega \subset \X$ be open, and let $u \in \W^{1,1}_{loc}(\Omega)$ be a solution of the minimal surface equation.
Let $\phi \in \Lip_c(\G(u))$, then ${\lip}_g(\phi) \cdot {\lip}_g(u_g) \in \sL^1(\G(u))$, and
\[
\int_{\G(u)} \lip_g(\phi) \cdot \lip_g(u_g) \, d\m_g=0.
\]
\end{customthm}

Finally, we outline the main ideas in the proofs of Theorems \ref{CT5} and \ref{CT6}. We will assume for simplicity that $u$ is locally Lipschitz on $\X$ and we will denote by $\nabla$ the relaxed gradient on $\X$ (for Lipschitz functions this coincides with the local Lipschitz constant $\lip$ thanks to \cite{Cheeger}). Both theorems can be easily obtained if one can show that given $\phi \in \Lip(\X)$ and the projection on the graph $i_g:\X \to \G(u)$, then the local Lipschitz constant of $\phi \circ i_g^{-1}$ (i.e. $\phi$ seen as a function on the graph $\G(u)$) satisfies
\[
    \lip_g(\phi \circ i_g^{-1})^2 \circ i_g=|\nabla \phi|^2-\frac{(\nabla \phi \cdot \nabla u)^2}{1+|\nabla u|^2}
    \quad \m \text{-a.e. on } \X.
\]
To obtain such an identity (which is what we would get on a Riemannian manifold) we use a technical blow-up argument so that we reduce the problem to the Euclidean case. \par 
To outline such argument we first need to recall a blow-up property of Lipschitz functions on metric measure spaces due to Cheeger (see \cite{Cheeger}).
Let $(\ssf{Y},\sd_y,\m_y)$ be a space with a locally doubling measure which supports a local Poincaré inequality (i.e. a $PI$ space) such that for $\m_y$-a.e. $y \in \ssf{Y}$ the tangent space of $\ssf{Y}$ at $y$ is a Euclidean space $\bb{R}^k$ (see Definition \ref{D4}). Given a point $y \in \ssf{Y}$ of the previous type and a sequence $\epsilon_n \to 0$, passing to a subsequence, we denote by $\psi_n:\bar{B}^{\bb{R}^k}_n(0)  \to (Y,\epsilon_n^{-1}\sd_y,y)$ the Gromov Hausdorff maps realizing the blow-up (see Definition \ref{D3}).
Thanks to \cite{Cheeger} we have that for every $\theta \in \Lip(\ssf{Y})$, for $\m_y$-a.e. $y \in \ssf{Y}$ there exists a linear function $\theta^{\infty}:\bb{R}^k \to \bb{R}$, called the blow-up of $\theta$ at $y$, such that
\[
\|\epsilon_n^{-1}(\theta \circ \psi_n-\theta(y))-\theta^{\infty}\|_{\sL^{\infty}(B^{\bb{R}^k}_1(0))} \to 0.
\]
Moreover, we have that the relaxed gradient $\nabla_y \theta (y)$ of $\theta$ in $y$ (which coincides with the local Lipschitz constant $\lip(\theta)(y)$) and the local Lipschitz constant $\lip(\theta^{\infty})(0)$ of $\theta^\infty$ in $0$ coincide.
\par
In particular, this blow-up property of Lipschitz functions holds for $\X$ and $\G(u)$ as these are both $PI$ spaces (since $u$ is now assumed to be Lipschitz). 
So we fix $\phi \in \Lip(\X)$ and a sequence $\epsilon_n \to 0$. We pick a point $x \in \X$ where the blow-up of $\X$ is realized (modulo passing to a subsequence) by the
Gromov Hausdorff maps $\psi_n:\bar{B}^{\bb{R}^k}_n(0) \to (\X,\epsilon_n^{-1}\sd,x)$ and the functions $u$ and $\phi$ admit blow-ups $u^{\infty}$ and $\phi^{\infty}$ respectively.
Let then $j:\bb{R}^k \to \Graph(u^{\infty})$ be the projection on the graph. It turns out that the maps
\[
    \psi'_n :=i_g \circ \psi_n \circ j^{-1}:j(B^{\bb{R}^k}_1(0)) \to (i_g(B^{\X}_{\epsilon_n}(x)), \epsilon_n^{-1}\sd_g)
\]
are Gromov Hausdorff maps realizing the blow-up of $\G(u)$ at $(x,u(x))$. In particular this implies that the blow-up of $\phi \circ i_g^{-1}$ on $\G(u)$ is $\phi^{\infty} \circ j^{-1}$. From this and the fact that blow-ups of functions preserve the Lipschitz constant, we obtain that
\begin{align*}
& \lip_g(\phi \circ i_g^{-1})^2 ( i_g(x))=\lip(\phi^{\infty} \circ j^{-1})^2  \\
& =\lip(\phi^{\infty})^2-\frac{(\lip(\phi^{\infty}) \cdot \lip(u^{\infty}))^2}{1+\lip (u^\infty)^2}= 
\Big( |\nabla \phi|^2-\frac{(\nabla \phi \cdot \nabla u)^2}{1+|\nabla u|^2} \Big)(x).
\end{align*}
If $u$ is not Lipschitz, the above strategy needs to be modified to ensure both existence of a blow-up of $u$ at points of $\X$ and existence of blow-ups of functions in $\Lip(\G(u))$ at points in $\G(u)$. 
Note in addition that in both cases we need the Lipschitz constant of the blow-up to coincide with the local Lipschitz constant of the initial function. 
These issues are addressed by using a geometric property of perimeter minimizers, i.e. the fact that these admit tangent balls to their boundary (see \cite{Weak} for the proof of this property in the non-collapsed case). \par
Finally, let us mention that Theorem \ref{CT1} might hold in the more general framework of doubling spaces supporting a $(1,1)$-Poincaré inequality (as this is the case in the Riemannian setting, see \cite{Ding}). Our strategy though, fails in this more general context. The main issues are the potential lack of regularity of tangent spaces to the ambient space and to the perimeter minimizing set, the non-linearity of the Laplacian and the possible failure of Theorem \ref{CT4} (due to the absence of a formula for the area of a graph in this more general setting).
\par 
The note is organized as follows: Section \ref{S4} contains preliminaries, Section \ref{S5} contains the proof of Theorem \ref{CT4}, Section \ref{S2} contains the proofs of Theorems \ref{CT5} and \ref{CT6}, while the final section contains the proofs of Theorems \ref{CT1}, \ref{CT2} and \ref{CT3}.

\medskip

\textbf{Acknowledgements.}
I would like to thank Prof. Andrea Mondino for the many discussions we had about the content of this note and for his valuable advice and guidance throughout this project. \\
The author is funded by a scholarship granted by the Mathematical Institute of the University of Oxford.

\nomenclature[01]{\(\Lip(\Omega)\)}{Lipschitz functions on $\Omega$. See beginning of Section \ref{S4}.}
\nomenclature[02]{\(\Lip_{loc}(\Omega)\)}{Locally Lipschitz functions on $\Omega$. See beginning of Section \ref{S4}.}
\nomenclature[03]{\(\Lip_c(\Omega)\)}{Compactly supported Lipschitz functions on $\Omega$. See beginning of Section \ref{S4}.}
\nomenclature[04]{\(\lip(f)\)}{Local Lipschitz constant of $f$. Defined in \eqref{E66}.}
\nomenclature[05]{\(\sL(f)\)}{Lipschitz constant of $f$. Defined in \eqref{E66}.}
\nomenclature[06]{\(\sd_e\)}{Euclidean distance.}
\nomenclature[07]{\(\Epi(f)\)}{Epigraph of $f$. Defined in \eqref{E70}.}
\nomenclature[08]{\(\Epi'(f)\)}{Hypograph of $f$. Defined in \eqref{E70}.}
\nomenclature[09]{\(W_f\)}{Abbreviation for $\sqrt{1+\lvert \nabla f \rvert^2}$.}
\nomenclature[10]{\(\sd_\times\)}{Product distance. See beginning of Section \ref{S5}.}
\nomenclature[11]{\(\m_\times\)}{Product measure. See beginning of Section \ref{S5}.}
\nomenclature[12]{\(\aH^{h}\)}{Codimension $1$ spherical Hausdorff measure. See Definition \ref{D7}.}
\nomenclature[13]{\((\epsilon_n,\psi_n)\)}{Blow-up of a metric space. See Definition \ref{D3}.}
\nomenclature[14]{\((\epsilon_n,\psi_n,f^\infty)\)}{Blow-up of a function. See Definition \ref{D8}.}
\nomenclature[15]{\(\G(u)\)}{Boundary of the open representative of $\Epi(u)$, where $u$ is a solution of the minimal surface equation. See Definition \ref{D5}.}
\nomenclature[16]{\(\sd_g\)}{Product distance restricted to $\G(u) \times \G(u)$. See Definition \ref{D5}.}
\nomenclature[17]{\(\m_g\)}{Perimeter measure induced by $\Epi(u)$ on the Borel subsets of $\G(u)$. See Definition \ref{D5}.}
\nomenclature[18]{\(\lip_g(f)\)}{Local Lipschitz constant of a function $f:\G(u) \to \bb{R}$. See beginning of Section \ref{S9}.}
\nomenclature[19]{\(i_g\)}{Projection on $\G(u)$. Defined in \eqref{E71}.}
\nomenclature[20]{\(u_g\)}{Height function on $\G(u)$. Defined in \eqref{E72}.}
\nomenclature[21]{\(\rho_{(x,t)}\)}{$l_1$ distance in $\X \times \bb{R}$ from $(x,t)$. See Definition \ref{D6}.}
\nomenclature[22]{\(D_{\bar{x},r}\)}{Ball centered in $\bar{x}$ w.r.t. $\rho_{\bar{x}}$ in $\X \times \bb{R}$. See Definition \ref{D6}.}
\nomenclature[23]{\(A_{\bar{x},r}\)}{Projection of $\G(u) \cap D_{\bar{x},r}$ on $\X$. See Definition \ref{D6}.}
\nomenclature[24]{\(E_{\bar{x},t}\)}{Superlevel set of $u_g$. See Definition \ref{D6}.}
\nomenclature[25]{\(E'_{\bar{x},t}\)}{Sublevel set of $u_g$. See Definition \ref{D6}.}
\nomenclature[26]{\(\Osc_{x,r}(f)\)}{Oscillation of $f$. Defined in \eqref{E73}.}

\printnomenclature

\section{Preliminaries} \label{S4}
Throughout the note will work on metric measure spaces $\mms{X}{d}{m}$, where $(\ssf{X},\ssf{d})$ is a separable complete length space where balls are precompact and $\m$ is a non-negative Borel measure on $\ssf{X}$ which is finite on bounded sets and whose support is the whole $\ssf{X}$. Given an open set $\Omega \subset \ssf{X}$ we denote by $\Lip(\Omega)$, $\Lip_{loc}(\Omega)$ and $\Lip_c(\Omega)$ respectively Lipschitz functions, locally Lipschitz and Lipschitz functions with compact support in $\Omega$. If $f \in \Lip_{loc}(\Omega)$ and $x \in \Omega$ we define
\begin{equation} \label{E66}
\lip(f)(x):= \limsup_{y \to x} \frac{|f(x)-f(y)|}{\ssf{d}(x,y)} \quad \text{and} \quad
\ssf{L}(f):= \sup_{x,y \in \Omega} \frac{|f(x)-f(y)|}{\ssf{d}(x,y)}.
\end{equation}
Given a closed interval $I \subset \bb{R}$, we say that a rectifiable curve $\gamma: I \to \X$ is a geodesic if its length coincides with the distance between its endpoints. Unless otherwise specified we assume that geodesics have constant unit speed. Throughout the note $\sd_e$ will be the Euclidean distance in any dimension.
In many proofs we will use the notation $c_1, \, c_2,$ etc. for constants that are independent of the other quantities appearing in the statement that we are proving.

\subsection{Sobolev spaces and Gromov Hausdorff convergence}
We now recall some basic notions about Sobolev spaces in the setting of metric measure spaces, the main references being \cite{Cheeger}, \cite{AGS}, \cite{AGS2} and \cite{Giglimem}.

\begin{definition}
    Let $\mms{X}{d}{m}$ be a metric measure space, $\Omega \subset \ssf{X}$ an open set and let $p>1$. A function $f \in \sL^p(\Omega)$ is said to be in the Sobolev space $\W^{1,p}(\Omega)$ if there exists a sequence of locally Lipschitz functions $\{f_i\}_{i \in \bb{N}} \subset \Lip_{loc}(\Omega)$ converging to $f$ in $\sL^p(\Omega)$ such that
    \[
    \limsup_{i \to + \infty} \int_\Omega \lip(f_i)^p \, d \m < + \infty.
    \]
     A function $f \in \sL^p_{loc}(\Omega)$ is said to be in the Sobolev space $\W^{1,p}_{loc}(\Omega)$ if for every $\eta \in \Lip_c(\Omega)$ we have $f \eta \in \W^{1,p}(\Omega)$.
\end{definition}

For any $f \in \W^{1,p}(\Omega)$ one can define an object $|\nabla f|$ (a priori depending on $p$, but independent of the exponent in the spaces that we will work on) such that for every open set $A \subset \Omega$ we have
    \[
    \int_A |\nabla f|^p \, d \m =
    \inf \Big\{ \liminf_{n \to + \infty} \int_A \lip(f_n)^p \, d \m \Big| (f_n)_n \subset \sL^p(A) \cap \Lip_{loc}(A), \|f_n-f\|_{\sL^p(A)} \to 0 \Big\}.
    \]
    The quantity in the previous expression will be called $p$-Cheeger energy and denoted by $\Ch_p(f)$ while $|\nabla f|$ will be called relaxed gradient. Later we will often write $\Ch$ in place of $\Ch_2$ for simplicity of notation.
    We define $\|f\|_{\W^{1,p}(\Omega)}:=\|f\|_{\sL^p(\Omega)}+\Ch_p(f)$. One can check that with this norm the space $\W^{1,p}(\Omega)$ is Banach.
We now introduce functions of bounded variation following \cite{Mir} (see also \cite{ADM}).

\begin{definition}
    Let $\mms{X}{d}{m}$ be a metric measure space and let $\Omega \subset \ssf{X}$ an open set. A function $f \in \sL^1(\Omega)$ is said to be of bounded variation if there exists a sequence of of locally Lipschitz functions $\{f_i\}_{i \in \bb{N}} \subset \Lip_{loc}(\Omega)$ converging to $f$ in $\sL^1(\Omega)$ such that
    \[
    \limsup_{i \to + \infty} \int_\Omega \lip(f_i) \, d \m < + \infty.
    \]
    The space of such functions is denoted $\BV(\Omega)$. A function $f \in \sL^1_{loc}(\Omega)$ is said to be in $\BV_{loc}(\Omega)$ if for every $\eta \in \Lip_c(\Omega)$ we have $f \eta \in \BV(\Omega)$.
\end{definition}

    For any $f \in \BV(\Omega)$ and any open set $A \subset \Omega$ we define
    \[
    |Df|(A) =
    \inf \Big\{ \liminf_{n \to + \infty} \int_A \lip(f_n) \, d \m \Big| (f_n)_n \subset \sL^1(A) \cap \Lip_{loc}(A), \|f_n-f\|_{\sL^1(A)} \to 0 \Big\}.
    \]
    One can check that the quantity in the previous expression is the restriction to the open subsets of $\Omega$ of a finite measure. We define $\|f\|_{\BV(\Omega)}:=\|f\|_{\sL^1(\Omega)}+|Df|(\Omega)$. One can check that with this norm the space $\BV(\ssf{X})$ is Banach. A function $f$ belongs to $\W^{1,1}(\Omega)$ if $f \in \BV(\Omega)$ and $|Df| \ll \m$. In this case we denote by $|\nabla f|$ the density of $|Df|$ with respect to $\m$.

\begin{definition}
    Let $\mms{X}{d}{m}$ be a metric measure space and let $\Omega \subset \ssf{X}$ be an open set. For every $p \geq 1$ we denote by $\W^{1,p}_0(\Omega)$ the closure in $\W^{1,p}(\Omega)$ of $\Lip_c(\Omega)$.
\end{definition}


\begin{definition}
    A metric measure space $\mms{X}{d}{m}$ is infinitesimally Hilbertian if the space $\W^{1,2}(\ssf{X})$ is a Hilbert space.
\end{definition}

 If $\mms{X}{d}{m}$ is Infinitesimally Hilbertian and $\Omega \subset \ssf{X}$ an open set, for every $f,g \in \W^{1,1}_{loc}(\Omega)$ we define the measurable function $\nabla f \cdot \nabla g:\Omega \to \bb{R}$ by 
 \[
 \nabla f \cdot \nabla g:=\frac{|\nabla (f+g)|^2-|\nabla(f-g)|^2}{4}.
 \]


As a consequence of the infinitesimal Hilbertianity assumption, the previously defined product of gradients is bilinear in both entries.
We then define the Laplacian in the metric setting.

 \begin{definition}
 Let $\mms{X}{d}{m}$ be infinitesimally Hilbertian and let $\Omega \subset \ssf{X}$ be an open set.
    Let $f \in \W^{1,2}(\Omega)$. We say that $f \in D(\Delta,\Omega)$ if there exists a function $h \in \sL^2(\Omega)$ such that
    \[
    \int_\Omega gh \, d \m=-\int_\Omega \nabla g \cdot \nabla f \, d \m \quad \text{for any } g \in \W^{1,2}_0(\Omega).
    \]
    In this case we say that $\Delta f=h$ in $\Omega$.
\end{definition}

We also have the following more general definition.

\begin{definition}
Let $\mms{X}{d}{m}$ be infinitesimally Hilbertian and let $\Omega \subset \ssf{X}$ be an open set.
    Let $f \in \W^{1,1}_{loc}(\Omega)$ and let $\mu$ be a Radon measure on $\Omega$. We say that $\Delta f=\mu$ in $\Omega$ in distributional sense if
    \[
    \int_\Omega h \, d\mu=-\int_\Omega \nabla g \cdot \nabla f \, d \m \quad \text{for any } g \in \Lip_c(\Omega).
    \]
\end{definition}

Similarly we define what it means to be a solution of the minimal surface equation.

\begin{definition}
 Let $\mms{X}{d}{m}$ be infinitesimally Hilbertian and let $\Omega \subset \ssf{X}$ be an open set.
    We say that $f \in \W^{1,1}_{loc}(\Omega)$ solves the minimal surface equation on $\Omega$ if for every $\phi \in \Lip_c(\Omega)$ we have
    \[
    \int_\Omega \frac{\nabla f \cdot \nabla \phi}{\sqrt{1+|\nabla f|^2}} \, d\m=0.
    \]
\end{definition}

We now recall the main definitions concerning Gromov Hausdorff convergence, referring to \cite{YuBurago_1992}, \cite{Vil}, and \cite{GMS13} for an overview on the subject. 
\begin{definition}
    Let $(\X,\sd_x,x)$ and $(\ssf{Y},\sd_y,y)$ be pointed metric spaces and let $\delta >0$. we say that a map $f:\X \to \ssf{Y}$ is a $\delta$-$GH$ map if
    \begin{itemize}
        \item $f(x)=y.$ 
        \item $\sup_{a,b \in \X} \Big|\sd_x(a,b)-\sd_y(f(a),f(b)) \Big| \leq \delta$. 
        \item The image of $f$ is a $\delta$-net in $\ssf{Y}$.
    \end{itemize}
\end{definition}


\begin{definition} \label{D2}
    We say that a sequence of pointed metric spaces $(\X_n,\sd_n,x_n)$ converges to $(\X,\sd,x)$ in pointed Gromov Hausdorff convergence if for every $\delta,R>0$ there exists $N$ such that for every $n \geq N$ there exists a $\delta$-$GH$ map $f_n^{\delta,R}:(\bar{B}_R(x_n),\sd_n,x_n) \to (\bar{B}_R(x),\sd,x)$.
\end{definition}


The previous definition of pointed Gromov Hausdorff convergence differs slightly from the one given in \cite[p. 272]{YuBurago_1992}. Nevertheless, when working on length spaces, the two definitions coincide.

\begin{definition} \label{D4}
    Let $(\X,\sd)$ be a metric space and $x \in \X$. We denote by $\Tan_x(\X)$ the (possibly empty) collection of (isometry classes of) metric spaces that are pointed Gromov Hausdorff limits as $r \downarrow 0$ of the family $(\X,r^{-1}\sd,x)$.
\end{definition}




\subsection{General properties of RCD(K,N) spaces}

We now recall some properties of $\RCD(K,N)$ spaces, i.e. infinitesimally Hilbertian metric measure spaces with Ricci curvature bounded from below by $K \in \bb{R}$ and dimension bounded from above by $N \in [1,+\infty)$ in synthetic sense. \par
The Riemannian Curvature Dimension condition $\RCD(K, \infty)$ was introduced in \cite{AGS2} (see
also \cite{Giglimem,AGMR}) coupling the Curvature Dimension condition $\ssf{CD}(K,\infty)$, previously pioneered
in \cite{S1,S2} and independently in \cite{V}, with the infinitesimal Hilbertianity assumption. The finite dimensional counterpart $\RCD(K,N)$ is then obtained coupling the finite dimensional Curvature Dimension condition $\ssf{CD}(K,N)$ with the infinitesimal Hilbertianity assumption and was proposed in \cite{Giglimem}. For a complete introduction to the topic we refer to the survey \cite{Asurv} and the references therein. 
Let us mention that in the literature one can find also the (a priori weaker) $\RCD^*(K,N)$. It was proved in \cite{EKS, AMS15}, that $\RCD^*(K,N)$ is equivalent to the dimensional Bochner inequality. Moreover, \cite{CavallettiMilman} (see also \cite{Liglob}) proved that $\RCD^*(K,N)$ and $\RCD(K,N)$ coincide.
We now recall the properties that we will use later on in the note.

The $\RCD(K,N)$ condition implies that the measure is locally doubling (see \cite{S1}) and the validity of a Poincaré inequality (see \cite{PI}). In particular, if $f$ is a locally Lipschitz function on a $\RCD(K,N)$ space, its relaxed gradient coincides with its local Lipschitz constant $\lip(f)$ (see \cite[Theorem $12.5.1$]{HKST} after \cite{Cheeger}).






The next theorem follows from \cite{MondinoNaber,BScon}.
\begin{thm} \label{constancy}
    Let $(\X,\sd,\m)$ be an $\RCD(K,N)$ space. There exists $k \in \bb{N} \cap [1,N]$, called essential dimension of $\X$, such that for $\m$-a.e. $x \in \X$ we have $\Tan_x(\X)=\{(\bb{R}^k,\sd_e)\}$. Any such point will be called a regular point for $\X$.
\end{thm}




We now recall some properties of the heat flow in the $\RCD$ setting, referring to \cite{AGMR,AGS2} for the proofs of these results.
Given an $\RCD(K,N)$ space $(\X,\sd,\m)$, the heat flow $P_t:\sL^2(\X) \to \sL^2(\X)$ is the $\sL^2(\X)$-gradient flow of the Cheeger energy $\Ch$.
It turns out that one can obtain a stochastically complete heat kernel $p_t:\X \times \X \to [0,+\infty)$, so that
the definition of $P_t(f)$ can be then extended to $\sL^{\infty}$ functions by setting
\[
P_t(f)(x):=\int_\X f(y) p_t(x,y) \, d\m(y).
\]
The heat flow has good approximation properties, in particular if $f \in \W^{1,2}(\X)$, then $P_t (f) \to f$ in $\W^{1,2}(\X)$; while if $f \in \sL^{\infty}(\X)$, then $P_t f \in \Lip(\X)$ for every $t>0$.

The next proposition follows combining the contractivity estimates for the heat flow of \cite{GBH14} with a standard lower semicontinuity argument.


\begin{proposition} \label{P19}
    Let $(\X,\sd,\m)$ be an $\RCD(K,N)$ space, let $\Omega \subset \X$ be an open set and let $f \in \BV(\X)$. If $|Df|(\partial \Omega)=0$, then
    $
    \lim_{t \to 0} |DP_t(f)|(\Omega)=|Df|(\Omega).
    $
\end{proposition}

\subsection{Sets of finite perimeter and minimal sets}

\begin{definition}
    Let $(\X,\sd,\m)$ be an $\RCD(K,N)$ space. Let $E \subset \X$. We say that $E$ has locally finite perimeter if $1_E \in \BV_{loc}(\X)$. For every Borel subset $B \subset \X$. We denote $|D1_E|(B)$ by $P(E,B)$.
\end{definition}
When considering the perimeter as a measure under an integral sign we will use the notation $\Per(E,\cdot)$ instead of $P(E,\cdot)$.

The next result is the relative isoperimetric inequality on $\RCD(0,N)$ spaces.

\begin{proposition} \label{P9}
    Let $(\X,\sd,\m)$ be an $\RCD(0,N)$ space. There exists $C(N)>0$ such that if $E \subset \X$ has finite perimeter, $r>0$ and $x \in \X$, then
    \[
     \m(B_r(x))^{\frac{1}{N}} \min \{\m(E \cap B_r(x)),\m(B_r(x) \setminus E)\}^{\frac{N}{N-1}}  \leq CrP(E,B_r(x)).
    \]
\end{proposition}

We now turn our attention to minimal sets.

\begin{definition}
Let $(\X,\sd,\m)$ be an $\RCD(K,N)$ space.
Let $\Omega \subset \X$ be an open set.
    Let $E \subset \Omega$ be a set of locally finite perimeter. We say that $E$ is perimeter minimizing in $\Omega$ if for every $x \in \Omega$, $r>0$ and $F \subset  \Omega$ such that $F \Delta E \subset \subset B_r(x) \cap \Omega$ we have that $P(E,B_r(x) \cap \Omega) \leq P(F,B_r(x) \cap \Omega)$. If we say that $E$ is perimeter minimizing we implicitly mean 
    that $\Omega=\X$.
    
Similarly, $E$ is locally perimeter minimizing in $\Omega$ if for every $x \in \Omega$ there exists $r>0$ such that for every $F \subset \Omega$ such that $F \Delta E \subset \subset B_r(x) \cap \Omega$ we have $P(E,B_r(x) \cap \Omega) \leq P(F,B_r(x) \cap \Omega)$. 
    If we say that $E$ is locally perimeter minimizing we implicitly mean that $\Omega=\X$.
\end{definition}

The next theorem comes from \cite[Theorem $4.2$ and Lemma $5.1$]{Dens}.

\begin{thm} \label{T7}
    Let $(\X,\sd,\m)$ be an $\RCD(K,N)$ space. There exist $C,\gamma_0>0$ depending only on $K$ and $N$ such that the following hold. If $E \subset \X$ is a set minimizing the perimeter in an open set $\Omega \subset \X$, then, up to modifying $E$ on an $\m$-negligible set, for any $x \in \partial E$ and $0<r<1$ such that $B_{2r}(x) \subset \Omega$ we have
    \[
    \frac{\m(E \cap B_r(x))}{\m(B_r(x))}>\gamma_0, \quad \frac{\m(B_r(x) \setminus E)}{\m(B_r(x))}>\gamma_0
    \]
    and
    \[
    \frac{\m(B_r(x))}{Cr} \leq P(E,B_r(x)) \leq \frac{C\m(B_r(x))}{r}.
    \]
\end{thm}
From the previous result one deduces that locally perimeter minimizing sets admit both a closed and an open representative, and these have the same boundary which in addition is $\m$-negligible. Whenever we consider the boundary of a locally perimeter minimizing set, we will implicitly be referring to the boundary of its closed (or open) representative. The next corollary follows easily from Theorem \ref{T7}.

\begin{corollary} \label{C11}
Let $(\X,\sd,\m)$ be an $\RCD(K,N)$ space. For every $R>0$ there exists $\gamma>0$ depending only on $K$, $N$ and $R$ such that the following happens.
    Let $E \subset \X$ be the closed representative of a set minimizing the perimeter in $\Omega \subset \X$. For every $x \in E$ and $0<r<R$ such that $B_{2r} \subset \Omega$ we have
    \[
    \frac{\m(E \cap B_r(x))}{\m(B_r(x))}
    \geq
    \gamma.
    \]
\end{corollary}



The next proposition can be found in \cite[Theorem $5.1$]{Tens} (see also \cite{Brena22}). Given $\Omega \subset \X$ and $f:\Omega \to \bb{R}$ we denote 
\begin{equation} \label{E70}
\Epi(f):=\{ (x,t) \in \Omega \times \bb{R} :t > f(x) \} \quad \text{and} \quad  \Epi'(f):=\{ (x,t) \in \Omega \times \bb{R} :t < f(x) \}.
\end{equation}

\begin{proposition} \label{P1}
Let $(\X,\sd,\m)$ be an $\RCD(K,N)$ space.
    Let $f \in \W^{1,1}_{loc}(\X)$. For every Borel set $B \subset \X$ we have
    \[
    P(\Epi(f),B \times \bb{R})=\int_B \sqrt{1+ | \nabla f|^2} \, d\m.
    \]
\end{proposition}

Given $f \in \W^{1,1}_{loc}(\X)$ we will use the notation $W_f:=\sqrt{1+|\nabla f|^2}$.

\begin{proposition} \label{P8}
    Let $(\X,\sd,\m)$ be an $\RCD(K,N)$ space. Let $\Omega \subset \X$ be open and let $u \in \W^{1,1}_{loc}(\Omega)$ be continuous at $x_0$ and such that $\Epi(u)$ is locally perimeter minimizing in $\Omega \times \bb{R}$, then there exists an open set $A \subset \Omega$ containing $x_0$ such that
    for every $\phi \in \W^{1,1}(A) \cap \sL^{\infty}(A)$ compactly supported in $A$ we have
    \begin{equation} \label{MSE}
        \int_\Omega \frac{\nabla u \cdot \nabla \phi}{W_u} \, d\m=0.
    \end{equation}
    \begin{proof}
    Since $\Epi(u)$ is locally perimeter minimizing in $\Omega \times \bb{R}$ there exists a cilynder $B_\delta(x_0) \times (u(x_0)-\epsilon,u(x_0)+\epsilon) \subset \Omega \times \bb{R}$ where $\Epi(u)$ minimizes the perimeter. Moreover, modulo decreasing $\delta$, we may assume that $|u-u(x_0)|< \epsilon/2$ in $B_\delta(x_0)$ by the continuity hypothesis on $u$. We then set $A:=B_\delta(x_0)$. Note that if $\phi$ is as in the statement, then the graph of the function $u+t\phi$ restricted to $A$, for $t \in \bb{R}$ small enough, will be contained in the cylinder $B_\delta(x_0) \times (u(x_0)-\epsilon,u(x_0)+\epsilon)$.
    This implies that the function $f:\bb{R} \to \bb{R}$ given by
        \[
        f(t):=P(\Epi(u+t\phi),A \times \bb{R})
        \]
        has a minimum in zero.
        By Proposition \ref{P1} we can write $f$ as
        \[
        f(t)=\int_A \sqrt{1+|\nabla(u+t\phi)|^2} \, d\m,
        \]
        so that by standard arguments involving Dominated Convergence Theorem (see \cite[Theorem 2.27]{Folland}) we have that $f$ is smooth and that 
        \[
        f'(0)=\int_\Omega \frac{\nabla u \cdot \nabla \phi}{W_u} \, d\m.
        \]
        Since $f$ is smooth and has a minimum in zero we deduce that $f'(0)=0$ as desired.
    \end{proof}
\end{proposition}

\begin{remark} \label{R1}
    The same argument of the previous proposition combined with a truncation argument shows that
    if $u \in \W^{1,1}_{loc}(\Omega)$ is bounded and such that $\Epi(u)$ is perimeter minimizing in $\Omega \times \bb{R}$, then for every compactly supported $\phi \in \W^{1,1}(\Omega)$ we have
    \[
        \int_\Omega \frac{\nabla u \cdot \nabla \phi}{W_u} \, d\m=0.
    \]
\end{remark}

\subsection{Existence of tangent balls to perimeter minimizers} \label{S7}

In this section we sketch the main steps to prove Theorem \ref{T2}, which will be of crucial importance later on. In \cite[Proposition $6.44$]{Weak} the same fact is proved assuming that the space is non collapsing. The proof in our setting is essentially the same, and for this reason some details are omitted. For simplicity we will work under the assumption that $E \subset \X$ is a perimeter minimizer, although Theorem \ref{T2} is stated for local perimeter minimizers in an open set.

The proof of the next result mimics the one of
\cite[Proposition $5.4$]{Weak}.

\begin{proposition} \label{P12}
    Let $(\X,\sd,\m)$ be an $\RCD(K,N)$ space.
    Let $E \subset \X$ be a perimeter minimizing set and consider the distance function $\sd_{\bar{E}}:\X \to \bb{R}$ from $\bar{E}$. Then $\Delta \sd_{\bar{E}}$ is a Radon measure in $\X$ and we have $\Delta \sd_{\bar{E}} =\Per(E,\cdot)+\Delta \sd_{\bar{E}} \mres (\X \setminus \bar{E})$ in distributional sense.
    \begin{proof}
    We do the proof assuming that $K=0$, the general case requiring only a slight modification.
        Thanks to \cite[Theorem $5.2$ and Proposition $6.1$]{BPS} for $\lambda^1$-a.e. $r >0$ and for every positive $\phi \in \Lip_c(\X)$, it holds 
        \begin{equation} \label{E41}
        \int_{\{\sd_{\bar{E}}>r \}}\phi \, d \Delta \sd_{\bar{E}}+\int_{\{ \sd_{\bar{E}}>r \}}\nabla \phi \cdot \nabla \sd_{\bar{E}} \, d\m=-\int_\X\phi \, d\Per(\{\sd_{\bar{E}}>r\}).
        \end{equation}
        Combining \cite[Lemma $2.41$]{Weak} and the Coarea Formula,
        we obtain that there exists a Radon measure $\mu$ which is the weak limit (in duality with compactly supported continuous functions) of the Radon measures $\Per(\{\sd_{\bar{E}}>r_i\})$ for some sequence $r_i \downarrow 0$.
        We claim that $\mu \ll \Per(E,\cdot)$. 
        It is clear that $\mu$ is supported on $\partial E$, so it is sufficient to prove that there exists a constant $c$ such that for every $x \in \partial E$, considering a sequence $t_k$ that goes to zero, we have
        \[
        \limsup_{t_k \to 0} \frac{\mu(B_{t_k}(x))}{\Per(E,B_{t_k}(x))} \leq c.
        \]
        So fix $t_k$ and note that
        \[
        \mu(B_{t_k}(x)) \leq \liminf_{i \to + \infty} \Per(\{\sd_{\bar{E}}>r_i\},B_{t_k}(x)) .
        \]
        Combining again \cite[Lemma $2.41$]{Weak} and the Coarea Formula, the r.h.s. is controlled by $CP(E,B_{2t_k}(x))$, which by Theorem \ref{T7} is in turn controlled by $C'P(E,B_{t_k}(x))$. Hence, we get
        \[
        \mu(B_{t_k}(x)) \leq C'P(E,B_{t_k}(x)),
        \]
        as desired. \par 
        With this in mind we pass to the limit in \eqref{E41} along the sequence $r_i$ and obtain
        \[
        \lim_{i \to + \infty} \int_{\{ \sd_{\bar{E}}>r_i \}}\phi \, d \Delta \sd_{\bar{E}}+
        \int_{\X}\nabla \phi \cdot \nabla \sd_{\bar{E}} \, d\m
        =-\int_\X\phi \, d\mu.
        \]
        The first summand in the l.h.s. of the previous equation is a linear functional on the continuous compactly supported functions. 
        Moreover it has negative sign because of \cite[Theorem $5.2$]{lapb} (we recall that we are assuming $K=0$ for simplicity), so it is represented by a negative Radon measure $\nu$.
        This implies, by the same equation, that $ \sd_{\bar{E}}$ has measure valued Laplacian and that $\Delta \sd_{\bar{E}} \mres (\X \setminus \bar{E})=\nu$, $\Delta \sd_{\bar{E}} \mres \partial E=\mu$ and $\Delta \sd_{\bar{E}} \mres E=0$. To conclude we only need to prove that $\mu=\Per(E)$. Since we know that $\mu \ll \Per(E)$ it is enough to show that for $\Per$-a.e. $x \in \partial E$ we have
        \[
        \lim_{t_k \to 0} \frac{\mu(B_{t_k}(x))}{\Per(E,B_{t_k}(x))} =1.
        \]
        This follows by the same blow-up argument of \cite[Proposition $5.4$ Step 2]{Weak}.
    \end{proof}
\end{proposition}

The next result will play the role of
\cite[Lemma $6.12$]{Weak}.

\begin{proposition} \label{P30}
Let $(\X,\sd,\m)$ be an $\RCD(K,N)$ space.
    Let $E \subset \X$ be a perimeter minimizing set and consider the distance function $\sd_{\bar{E}}:\X \setminus \bar{E} \to \bb{R}$. Then for every non-negative $\phi \in \Lip_c(\X)$ and $\lambda^1$-a.e. $r \in \bb{R}$, it holds
    \[
    \int_\X \phi \, d \Per(\{\sd_{\bar{E}} >r\})-\int_\X \phi \, d \Per(E) \leq \int_{ \{0 < \sd_{\bar{E}} < r \}} \nabla \phi \cdot \nabla \sd_{\bar{E}}+\phi t_{K,N}(\sd_{\bar{E}}) \, d\m,
    \]
    where $t_{K,N} \in C^\infty([0,+\infty))$ and $t_{K,N}(0)=0$.
    \begin{proof}
        For every $r>0$ define $E_r:=\{\sd_{\bar{E}}<r\}$. Applying \cite[Theorem $5.2$ and Proposition $6.1$]{BPS} with $A:=E_r$ and $g=\sd_{\bar{E}}$ (this can be done thanks to Proposition \ref{P12}) we get that for every non-negative $\phi \in \Lip_{loc}(\X)$ we have
        \begin{align*}
        \int_{\partial E} \phi \, d\Delta \sd_{\bar{E}}+\int_{E_r \setminus \bar{E}} \phi \,d\Delta \sd_{\bar{E}}+\int_{\{0<\sd_{\bar{E}}<r\}} \nabla \phi \cdot \nabla \sd_{\bar{E}} \, d\m
        = \int_{\X}\phi \, d\Per(E_r).
        \end{align*}
        By \cite[Theorem 5.2]{lapb}, there exists $t_{K,N} \in C^\infty([0,+\infty))$ with $t_{K,N}(0)=0$, such that $\Delta \sd_{\bar{E}} \mres (\X \setminus \bar{E}) \leq t_{K,N}(\sd_{\bar{E}}) \m$.
        The statement then follows by Proposition \ref{P12}. 
    \end{proof}
\end{proposition}

The proof of the next result mimics the one of
\cite[Proposition $6.14$]{Weak}.

\begin{proposition} \label{P32}
Let $(\X,\sd,\m)$ be an $\RCD(K,N)$ space.
    Let $E \subset \X$ be a perimeter minimizing set and consider the measure $\mu_\epsilon:=\epsilon^{-1}\m \mres \{0<\sd_{\bar{E}}<\epsilon\}$. Then $\mu_\epsilon$ converges weakly to $\Per(E)$ as $\epsilon \to 0$.
    \begin{proof}
        Fix a sequence $\epsilon_k \downarrow 0$. By the same argument of \cite[Proposition 6.14, first lines]{Weak}, there exists a limit measure $\mu$ such that $\mu_{\epsilon_k} \to \mu$ in weak sense and $\mu \geq \Per(E)$. To prove the reverse inequality, we mimic the proof of \cite[Proposition 6.14]{Weak}, where Proposition \ref{P30} replaces \cite[Lemma 6.12]{Weak}.
        To this aim fix a non-negative function $\phi \in \Lip_c(\X)$. By Proposition \ref{P30}, for $\lambda^1$-a.e. $r \in \bb{R}$, it holds
    \[
    \int_\X \phi \, d \Per(\{\sd_{\bar{E}} >r\})-\int_\X \phi \, d \Per(E) \leq \int_{\{0 < \sd_{\bar{E}} < r\}} \nabla \phi \cdot \nabla \sd_{\bar{E}} +\phi t_{K,N}(\sd_{\bar{E}}) \, d\m.
    \]
    Hence, using the Coarea Formula and the previous inequality,
    \begin{align*}
        \int_\X & \phi \, d\mu=\lim_{s \to 0} \fint_0^s
        \int_\X \phi \, d\Per(\{\sd_{\bar{E}} >r\}) \, dr \\
        & \leq \int_\X \phi \, d \Per(E)+\limsup_{s \to 0} \fint_0^s \int_{ \{0 < \sd_{\bar{E}} < r\}} \nabla \phi \cdot \nabla \sd_{\bar{E}} +\phi t_{K,N}(\sd_{\bar{E}}) \, d\m \, dr=\int_\X \phi \, d \Per(E),
    \end{align*}
    concluding the proof.
    \end{proof}
\end{proposition}

The proof of the next result mimics the one of
\cite[Proposition $6.15$]{Weak}.

\begin{proposition} \label{P33}
Let $(\X,\sd,\m)$ be an $\RCD(K,N)$ space.
     Let $E \subset \X$ be a perimeter minimizing set. Let $h>0$, let $\Gamma \subset \{\sd_{E}=h\}$ be compact, and set $\Gamma_E:=\{y \in \partial E: \exists x \in \Gamma : \sd(y,x)=h\}$. Then, setting $G:=\{x : \sd_E(x)+\sd_\Gamma(x)=h\}$, it holds $P(\{\sd_E < h\},\Gamma) \leq P(E,\Gamma_E)+ \int_G t_{K,N}(\sd_E) \, d\m$, where $t_{K,N} \in C^\infty([0,+\infty))$ and $t_{K,N}(0)=0$.
     \begin{proof}
         Let $\epsilon>0$ and set $\phi_\epsilon:=\epsilon^{-1}(h+\epsilon-\sd_E-\sd_\Gamma) \vee 0$. For any $\delta \in (0,h)$ set $g_\delta$ to be the function such that $g_\delta(0)=g_\delta'(0)$ and $g_\delta''=\delta^{-1}(\chi_{[0,\delta]}-\chi_{[h-\delta,h]})$. By \cite[Theorem 5.2]{lapb} and the chain rule, it holds $\Delta g_\delta(\sd_E) \leq g''_\delta(\sd_E)+g_\delta'(\sd_E) t_{K,N}(\sd_E)$. 

         Consider an open neighbourhood $F$ of $G$. By \cite[Theorem $5.2$]{BPS} (taking into account the simplifications explained in \cite[lines after (6.15)]{Weak}), one deduces that, for $\epsilon>0$ small enough,
         \[
         \int_F g_\delta''(\sd_E) \phi_\epsilon \, d\m \geq -\int_F g_\delta'(\sd_E) t_{K,N}(\sd_E) \phi_\epsilon \, d\m.
         \]
         By definition of $g_\delta$, it holds
         \[
         \delta^{-1} \int_F \chi_{[0,\delta]}(\sd_E) \phi_\epsilon \, d\m \geq \delta^{-1} \int_F \chi_{[h-\delta,h]}(\sd_E) \phi_\epsilon \, d\m-\int_F g_\delta'(\sd_E) t_{K,N}(\sd_E) \phi_\epsilon \, d\m.
         \]
         Passing to the limit as $\delta \to 0$, using Proposition \ref{P32}, the lower semicontinuity of total variations and the fact that $g_\delta' \to 1$ as $\delta \to 0$, we deduce
         \[
         \int_F \phi_\epsilon \, d\Per(E) \geq \int_F \phi_\epsilon \, d \Per (\{\sd_E < h\})-\int_F t_{K,N}(\sd_E) \phi_\epsilon \, d\m.
         \]
         Since $\phi_\epsilon \to \chi_G$ as $\epsilon \to 0$, passing to the limit in the previous inequality we conclude.
     \end{proof}
\end{proposition}


The proof of the next result mimics the one of
\cite[Proposition $6.44$]{Weak}.

\begin{thm} \label{T2}
Let $(\X,\sd,\m)$ be an $\RCD(K,N)$ space.
    Let $\Omega \subset \X$ be an open set and let $E \subset \Omega$ be a set locally minimizing the perimeter in $\Omega$. For $\Per$-a.e. $x \in \partial E \cap \Omega$ there exists balls $B_1 \subset E$ and $B_2 \subset \Omega \setminus E$ such that $\partial B_1 \cap \partial B_2=\{x\}$. These balls are called tangent balls to $E$ at $x$.
    \begin{proof}
    For simplicity we assume that $E \subset \X$ is perimeter minimzing.
        Let $x \in \partial E$. By \cite[Corollary 6.42]{Weak} (the proof works exactly in the same way in our setting), there exist $c(K,N),r_0(K,N)>0$ such that for every $0<r<r_0$ there exists $1/2<\rho<1$ such that 
        \begin{equation} \label{E47}
        P(B_\rho(x),\{0<\sd_E<r\}) \leq crP(E,B_2(x)), \quad P(\{0<\sd_E<r\},B_\rho(x))=0.
        \end{equation}
        Set now $\Gamma_r:=\{\sd_E=r\} \cap B_\rho(x)$ and let $\Gamma_{r,E}:=\{y \in \partial E: \exists z \in \Gamma_r :\sd(y,z)=r\}$.
        Since $E$ is perimeter minimizing, it holds
        \begin{align*}
        P(E,B_{1/2}(x)) & \leq P(E \cup (\{\sd_E < r\} \cap B_\rho(x)),B_{1/2}(x)), 
        \end{align*}
        so that combining with \eqref{E47} and Proposition \ref{P33}, it holds
        \begin{align*}
        P(E,B_\rho(x)) & \leq P(\{\sd_E<r\},B_\rho(x))+P(B_\rho(x),\{0<\sd_E<r\}) \\
        &  \leq P(E,\Gamma_{r,E})+crP(E,B_2(x))+o(r).
        \end{align*}
        Note that $\Gamma_{r,E}$ is contained in $A_r^e:=\{y \in \partial E: \exists \text{ a ball } B_r(z) \subset \X \setminus E:y \in \partial B_r(z)\}$. Hence
        \[
        P(E,B_{1/2}(x) \setminus A^e_r) \leq  P(E,B_\rho(x) \setminus A^e_r) \leq crP(E,B_2(x))+o(r).
        \]
        In particular, $\Per(E)$-a.e. $y \in \partial E$ admits a tangent ball from outside. By a similar argument, $\Per(E)$-a.e. $y \in \partial E$ admits also a tangent ball from the inside.
    \end{proof}
\end{thm}


\section{Minimal surface equation and perimeter minimizers} \label{S5}

We will denote by $\sd_\times$ and $\m_\times$ respectively the product distance and the product measure in $\X \times \bb{R}$. 
We recall that given $\Omega \subset \X$ and $f:\Omega \to \bb{R}$ we denote by $\Epi(f)$ the set $\{ (x,t) \in \Omega \times \bb{R} :t > f(x) \}$ and by $\Epi'(f)$ the set $\{ (x,t) \in \Omega \times \bb{R} :t < f(x) \}$.
The goal of this section is to prove the following theorem (which coincides with Theorem \ref{CT4} in the Introduction when $\Omega=\X$). 

\begin{thm} \label{T11}
    Let $(\X,\sd,\m)$ be an $\RCD(K,N)$ space, let $\Omega \subset \X$ be open and let $u \in \W^{1,1}_{loc}(\Omega)$. The following conditions are equivalent.
    \begin{enumerate}
        \item $\Epi(u)$ is locally perimeter minimizing in $\Omega \times \bb{R}$.
        \item $\Epi(u)$ is perimeter minimizing in $\Omega \times \bb{R}$.
        \item $u$ solves the minimal surface equation on $\Omega$.
    \end{enumerate}
\end{thm}

In Subsection \ref{SS1} we show that $1$ implies $3$, while in Subsection \ref{SS2} we will show that $3$ implies $2$.


\subsection{Locally perimeter minimizing implies MSE} \label{SS1}

The proof of the next result is inspired from \cite[Theorem $14.10$]{Giusti}. We recall that whenever we refer to a locally perimeter minimizing set we implicitly mean its open representative. Moreover
whenever we refer to the pointwise behavior of $u$ we mean its precise representative defined by
        \begin{equation} \label{E45}
        u(x):=\limsup_{r \to 0}\fint_{B_r(x)}u \, d\m.
        \end{equation}

\begin{proposition} \label{P14}
    Let $(\X,\sd,\m)$ be an $\RCD(K,N)$ space and let $\Omega \subset \X$ be open. Let $u \in \W^{1,1}_{loc}(\Omega)$ be such that $\Epi(u)$ is locally perimeter minimizing in $\Omega \times \bb{R}$, then $u$ is continuous at its Lebesgue points.
    \begin{proof}
        Let $x \in \Omega$ be a Lebesgue point for $u$ and suppose that $u(x) = 0$. Modulo a vertical translation we always trace back to this case. Since $\Epi(u)$ minimizes the perimeter locally, there exists an open set of the form $B_\epsilon(x) \times (-\epsilon,\epsilon) \subset \Omega \times \bb{R}$ where $\Epi(u)$ is perimeter minimizing. Now consider $r>0$ such that $\epsilon >4r$ and $y \in B_r(x)$. Suppose for now that $u(y) \geq 0$. 
        \par For every $i \in \bb{N}$ such that $2ir \in [0,\min\{u(y),\epsilon /4 \}]$ we have that $(y,2ir)$ is in the closure of $\Epi'(u)$ and $B_{r}(y,2ir) \subset \X \times \bb{R}$ is contained in $B_{\epsilon/2} \times (-\epsilon /2,\epsilon/2)$.

        We denote by $c_1,c_2, \cdots$ positive constants depending on $K$ and $N$.
        Applying Corollary \ref{C11} in the product space $(\X \times \bb{R},\sd_\times,\m_\times)$, it holds
        \[
        \m_\times(B_{r}(y,2ir) \cap \Epi'(u)) \geq 
        c_1\m_\times(B_r(y,2ir)) \geq c_2r\m(B_r(y)).
        \]
        In particular, it holds
        \begin{align*}
        \int_{B_{2r}(x)} |u| \, d\m & \geq 
        \sum_{i \in \bb{N} \cap [0,\min\{u(y)/2r,\epsilon/8r\}]} \m_\times(B_{r}(y,2ir) \cap \Epi'(u)) 
        \\
        & \geq c_3r(\min\{u(y)/2r,\epsilon/8r\}-2)\m(B_r(y)).
        \end{align*}
        By the doubling property,
        \[
        c_3r+\fint_{B_{2r}(x)} |u| \, d\m \geq
        c_5\min\{u(y),\epsilon\}.
        \]
        In particular, since $u(y)$ was assumed to be positive and $x$ is a Lebesgue point with $u(x)=0$, we deduce that
        $
         \sup_{B_r(x)}u \vee 0 \to 0
        $
        as $r \downarrow 0$.
        By an analogous argument, it holds
        $
         \inf_{B_r(x)}u \wedge 0 \to 0
        $
        as well,
        proving the continuity of $u$ in $x$.
    \end{proof}
\end{proposition}

\begin{definition} \label{D7}
    Given an $\RCD(K,N)$ space $(\X,\sd,\m)$ we define the codimension $1$ spherical Hausdorff measure to be the measure obtained with the Carathéodory construction using coverings made by balls and gauge function 
    \[
    B_r(x) \mapsto \m(B_r(x))/r.
    \]
    We denote such measure by $\aH^{h}$.
\end{definition}

The next proposition corresponds to \cite[Theorem $4.1$]{Lebes}.

\begin{proposition} \label{P15}
    Let $(\X,\sd,\m)$ be an $\RCD(K,N)$ space and let $\Omega \subset \X$ be open. Let $f \in \W^{1,1}_{loc}(\Omega)$. Then $\aH^{h}$-a.e. $x \in \Omega$ is a Lebesgue point of $f$.
\end{proposition}

The next proposition is an intermediate step to prove that $3$ implies $1$ in Theorem \ref{T11}.

\begin{proposition} \label{P16}
    Let $(\X,\sd,\m)$ be an $\RCD(K,N)$ space and let $\Omega \subset \X$ be open. Let $u \in \W^{1,1}_{loc}(\Omega)$ be such that $\Epi(u)$ is locally perimeter minimizing in $\Omega \times \bb{R}$. There exists an open set $A \subset \Omega$ with $\aH^{h}(\Omega \setminus A)=0$ such that, for every $\phi \in \Lip_c(A)$, it holds
    \[
    \int_\Omega \frac{\nabla u \cdot \nabla \phi}{W_u} \, d\m=0.
    \]
    \begin{proof}
        Combining Propositions \ref{P14} and \ref{P15} we get that $u$ is continuous $\aH^{h}$-almost everywhere. For every continuity point $x$ of $u$ consider the set $A_x \subset \Omega$ given by Proposition \ref{P8} and denote with $A$ the union of these sets. It is clear that $\aH^{h}(\Omega \setminus A)=0$. 
        Now let $\phi \in \Lip_c(A)$ and let $\{A_{x_i}\}_{i=1}^m$ be a finite cover of 
        the support of $\phi$. Let $\{\eta_i\}_{i=1}^m$ be Lipschitz functions such that their sum is equal to $1$ on the support of $\phi$, while each $\eta_i$ is compactly supported in $A_{x_i}$. 
        It is easy to check that such functions exist. 
        We get
        \[
         \int_\Omega \frac{\nabla u \cdot \nabla \phi}{W_u} \, d\m
         = \sum_{i=1}^m 
         \int_{A_i} \frac{\nabla u \cdot \nabla (\eta_i \phi)}{W_u} \, d\m=0.
        \]
    \end{proof}
\end{proposition}

\begin{proposition} \label{Pf}
    Let $(\X,\sd,\m)$ be an $\RCD(K,N)$ space and let $\Omega \subset \X$ be open. Let $u \in \W^{1,1}_{loc}(\Omega)$ be such that $\Epi(u)$ is locally perimeter minimizing in $\Omega \times \bb{R}$. For every $\phi \in \Lip_c(\Omega)$, it holds
    \[
    \int_\Omega \frac{\nabla u \cdot \nabla \phi}{W_u} \, d\m=0.
    \]
    \begin{proof}
        Let $A$ be the set given by Proposition \ref{P16} and call $C:=\Omega \setminus A$. We claim that we can construct a sequence $\{\eta_i\}_{i \in \bb{N}}$ of Lipschitz compactly supported functions in $\Omega$ that are equal to $1$ in a neighbourhood of $C \cap \supp(\phi)$ and such that $\|\eta_i\|_{\W^{1,1}(\Omega)} \to 0$ as $i$ goes to $+\infty$. \par
        Assume for the moment that the claim holds. In this case, it holds
        \begin{align*}
        \int_\Omega \frac{\nabla u \cdot \nabla \phi}{W_u} \, d\m &=\int_\Omega \frac{\nabla u \cdot \nabla (\eta_i\phi)}{W_u} \, d\m
        +\int_\Omega \frac{\nabla u \cdot \nabla ((1-\eta_i)\phi)}{W_u} \, d\m
        \\
        &=\int_\Omega \frac{\phi\nabla u \cdot \nabla \eta_i}{W_u} \, d\m + 
        \int_\Omega \frac{\eta_i \nabla u \cdot \nabla \phi}{W_u} \, d\m.
        \end{align*}
        Since $|\nabla u|/W_u \leq 1$, the last expression tends to zero as $i$ goes to $+\infty$ since $\|\eta_i\|_{\W^{1,1}(\Omega)} \to 0$. \par 
        So we only need to prove our previous claim. We have that $C \cap \supp(\phi)$ is compact and $\aH^{h}(C)=0$. So for every $\epsilon>0$ there exists a finite collection $\{B_{r_i^{\epsilon}}(x_i^\epsilon)\}_{i=1}^{m_\epsilon}$ of balls with radii in $(0,\epsilon)$ and centers in $C \cap \supp(\phi)$, whose union covers $C \cap \supp(\phi)$, and such that
        \[
        \sum_{i=1}^{m_\epsilon} \frac{\m(B_{r_i^{\epsilon}}(x_i^\epsilon))}{r_i^{\epsilon}} <\epsilon.
        \]
        For every such ball define $\eta_i^{\epsilon}:\X \to \bb{R}$ by
        \[
        \eta_i^{\epsilon}(x)=\Big(1-\frac{\sd(B_{r_i^{\epsilon}}(x_i^\epsilon),x)}{r_i^{\epsilon}} \Big) \vee 0.
        \]
        It is clear that, for some constant $c_1(K,N)>0$, it holds $\|\eta_i^{\epsilon}\|_{\sL^1(\X)} \leq \m(B_{2r_i^{\epsilon}}(x_i^\epsilon)) \leq c_1 \m(B_{r_i^{\epsilon}}(x_i^\epsilon))$, and
        \[
        \|\nabla \eta_i^{\epsilon}\|_{\sL^1(\X)} \leq \frac{ \m(B_{2r_i^{\epsilon}}(x_i^\epsilon))}{r_i^{\epsilon}} \leq c_1 \frac{ \m(B_{r_i^{\epsilon}}(x_i^\epsilon))}{r_i^{\epsilon}}.
        \]
        In particular defining $\eta_\epsilon:=(\sum_{i=1}^{m_\epsilon}\eta_i^{\epsilon}) \wedge 1$ we obtain a Lipschitz function compactly supported in $\Omega$ (if $\epsilon$ is small enough) which is equal to $1$ on a neighbourhood of $C \cap \supp(\phi)$ and with $\|\eta_\epsilon\|_{\W^{1,1}(\X)} \leq 2c_1 \epsilon$. This concludes the proof.
    \end{proof}
\end{proposition}

\begin{remark}
    The result of Proposition \ref{Pf} holds also in the setting of doubling spaces supporting a $(1,1)$ Poincaré inequality (without any curvature assumption). This follows taking into account that the density estimates for perimeter minimizers (that were used in Proposition \ref{P14}) hold in this weaker setting as well.
\end{remark}

\subsection{MSE implies globally perimeter minimizing} \label{SS2}

So far we have seen that, if $u \in \W^{1,1}_{loc}(\Omega)$ and its epigraph is locally perimeter minimizing in $\Omega \times \bb{R}$, then $u$ solves the minimal surface equation. In this section we show that if $u$ solves the aforementioned equation, then its epigraph minimizes the perimeter.

\begin{proposition}
    Let $(\X,\sd,\m)$ be an $\RCD(K,N)$ space and let $\Omega \subset \X$ be open and bounded. Let $u \in \W^{1,1}(\Omega)$ be a solution of the minimal surface equation. Then, for every $v \in \W^{1,1}(\Omega)$ such that $\{u \neq v\} \subset \subset \Omega$, it holds 
    \[
    \int_\Omega \sqrt{1+|\nabla u|^2} \, d\m
    \leq \int_\Omega \sqrt{1+|\nabla v|^2} \, d\m.
    \]
    \begin{proof}
        The function $f(t):=\int_A \sqrt{1+|\nabla(u+t(v-u))|^2} \, d\m$ belongs to $C^\infty(\bb{R})$ by \cite[Theorem 2.27]{Folland}. Since $u$ solves the minimal surface equation, it holds $f'(0)=0$. By taking second derivatives and using again \cite[Theorem 2.27]{Folland}, it is easy to check that $f$ convex, so that it has a minimum in $0$. 
    \end{proof}
\end{proposition}

The previous proposition implies that the theory of De Giorgi Classes can be applied to $u$. The next proposition follows by mimicking the proof in the Euclidean setting given for example in \cite[Theorem $3.9$]{Beck} (and it can be also obtained combining the results from \cite{HKLreg} and \cite{Tens}).

\begin{proposition}
     Let $(\X,\sd,\m)$ be an $\RCD(K,N)$ space and let $\Omega \subset \X$ be open. Let $u \in \W^{1,1}_{loc}(\Omega)$ be a solution of the minimal surface equation, then $u \in \sL^\infty_{loc}(\Omega)$.
\end{proposition}

The next proposition contains a technical approximation result that will be crucial for the remaining part of the section.

\begin{proposition} \label{P21}
     Let $(\X,\sd,\m)$ be an $\RCD(K,N)$ space. Let $\Omega' \subset \subset \Omega'' \subset \subset \Omega$ be open sets and let $f \in \W^{1,1}(\Omega)$ and $\phi \in \BV_c(\Omega')$ be both bounded. There exist functions $f_t \in \Lip(\Omega'')$ converging in $\sL^1(\Omega'')$ to $f+\phi$ such that 
    \begin{equation} \label{E52}
    \limsup_{t \to 0} P(\Epi(f_t),\Omega'' \times \bb{R}) = P(\Epi(f+\phi),\Omega'' \times \bb{R})
    \end{equation}
    and 
    \begin{equation} \label{E53}
    |D(f-f_t)|(\Omega''\setminus \Omega') \to 0.
    \end{equation}
    \begin{proof}
    Since the statement concerns precompact sets of $\Omega$, modulo using a cut off, we suppose that $f \in \W^{1,1}(\Omega)$ with compact support in $\Omega$. Moreover, modulo enlarging $\Omega''$, we may also suppose that $\m(\partial \Omega'')=0$. Whenever we refer to $f$ as a function on $\X$ we implicitly mean its extension to zero. 
    
    Let $B$ be an open set such that  $\supp(\phi) \subset \subset B \subset \subset \Omega'$.
        We claim that, for every $\epsilon >0$, there exists a function $f'\in \W^{1,1}(\Omega)$ with compact support, whose restriction to $\Omega \setminus B$ is Lipschitz, and such that $\|f'-f\|_{\sL^1(\Omega'')} \leq \epsilon$,
        \[
        |P(\Epi(f+\phi),\Omega'' \times \bb{R})-P(\Epi(f'+\phi),\Omega'' \times \bb{R})| < \epsilon,
        \]
        and
        \[
        |D(f-f')(\Omega'' \setminus \Omega')|< \epsilon.
        \]
        To prove the claim
        consider an open set $A$ such that
        \[
        \supp(\phi) \subset \subset A \subset \subset B, 
        \]
        let $\delta>0$, and let $v \in \Lip_c(\Omega)$ be such that $\|f-v\|_{\W^{1,1}(\Omega)}< \delta$. 
        Let then $\eta \in \Lip_c(B)$ be a positive function, taking value less than or equal to $1$ and identically equal to $1$ on $A$. 
        We define $f':=\eta f + (1-\eta) v$. This function trivially satisfies $f' \in \W^{1,1}(\Omega)$, it has compact support, its restriction belongs to $\Lip (\Omega \setminus B)$, $\|f'-f\|_{\sL^1(\Omega'')}+|D(f-f')(\Omega'' \setminus \Omega')|< c(\eta)\delta$, and
        \[
         |P(\Epi(f+\phi),\Omega'' \times \bb{R})-P(\Epi(f'+\phi),\Omega'' \times \bb{R})|
         \]
         \[
         =
         |P(\Epi(f),{^c}A \times \bb{R})-P(\Epi(f'),{^c}A \times \bb{R})|.
         \]
         By Proposition \ref{P1} and triangular inequality, this last quantity is controlled by
         \[
         |D(f'-f)|({^c}A) 
         \leq |D(1-\eta)(v-f)|(\Omega) 
         \leq c(\eta) \delta.
         \]
         In particular, choosing $\delta$ small enough we have proved the claim. 
        Taking this into account, it is sufficient to prove the statement of the theorem under the additional assumption that the restriction of $f$ to $\Omega \setminus B$ is Lipschitz (and we still have that $f \in \W^{1,1}(\Omega)$ and has compact support). \par 
        To this aim we define $f_t:=P_t(f+\phi)$ and we claim that these functions restricted to $\Omega''$ have the right properties.
        First of all the functions $f_t$ are Lipschitz by the $\sL^{\infty}$ to Lipschitz property of the heat flow. Moreover $f_t \to f+\phi$ in $\sL^1(\Omega'')$ since we have convergence in $\sL^2(\X)$. \par 
        We will now show that
        \[
    |D(f-f_t)|(\Omega''\setminus \Omega') \to 0.
    \]
    Let $\tau \in \Lip_c(\Omega')$ be positive, taking value less than or equal to $1$ and equal to $1$ on $B$.
        Note that 
        \[
        f_t=P_t((1-\tau)f)+P_t(\tau f+\phi),
        \]
        and since $(1-\tau)f \in \Lip(\X) \subset \W^{1,2}(\X)$ we have $P_t( (1-\tau)f) \to (1-\tau)f$ in $\W^{1,2}(\X)$ which implies convergence in $\W^{1,1}(\Omega'')$ since $\Omega''$ is bounded. 
        In particular, when we restrict these functions to $\Omega'' \setminus \Omega'$ we get
        \[
        |D(f-P_t((1-\tau)f))|(\Omega''\setminus \Omega') \to 0.
        \]
        So to conclude the proof of \eqref{E53} it is sufficient to note that, by Proposition \ref{P19}, it holds
        \[
        |D(P_t(\tau f+\phi))|(\Omega''\setminus {\Omega}') \to |D(\tau f+\phi)|(\Omega''\setminus {\Omega}')=0.
        \]
        Finally, condition \eqref{E52} follows from \cite[Theorem 1]{AreaFormula}.
    \end{proof}
\end{proposition}

\begin{proposition} \label{P22}
    Let $(\X,\sd,\m)$ be an $\RCD(K,N)$ space and let $\Omega \subset \X$ be open. Let $u \in \W^{1,1}_{loc}(\Omega)$ be a solution of the minimal surface equation. Then the epigraph of $u$ minimizes the perimeter among bounded competitors in $\BV_{loc}(\Omega)$ that coincide with $u$ out of a compact set.
    \begin{proof}
    Let $\Omega' \subset \subset \Omega'' \subset \subset \Omega$ be open sets and let $\epsilon>0$. Let $\phi \in \BV_c(\Omega')$ be bounded. 
        It will be sufficient to show that for a constant $C$ independent of $\epsilon$ we have
        \[
        P(\Epi(\phi+u),\Omega'' \times \bb{R}) \geq P(\Epi(u),\Omega'' \times \bb{R})-C\epsilon.
        \] 
         By Proposition \ref{P21}, there exists $f \in \Lip(\Omega'')$ such that
         \[
         \|f-u\|_{\W^{1,1}(\Omega'' \setminus \bar{\Omega}')}+\|f-(u+\phi)\|_{\sL^1(\Omega'')} < \epsilon
         \]
         and
         \[
         P(\Epi(f),\Omega'' \times \bb{R}) \leq P(\Epi(\phi+u),\Omega'' \times \bb{R})+\epsilon.
         \]
         Let $\eta \in \Lip_c(\Omega'')$ be positive, taking value less than or equal to $1$ and equal to $1$ on a neighbourhood of $\bar{\Omega}'$. 
         Note that $f \eta + (1-\eta)u \in \W^{1,1}_{loc}(\Omega)$ and differs from $u$ on a precompact set of $\Omega''$, so that
         \[
         P(\Epi(u),\Omega'' \times \bb{R}) \leq 
         P(\Epi(f \eta + (1-\eta)u),\Omega'' \times \bb{R}).
         \]
         At the same time
         \[
         |P(\Epi(f \eta + (1-\eta)u),\Omega'' \times \bb{R})-P(\Epi(f),\Omega'' \times \bb{R})|
         \leq |D(1-\eta)(f-u)|(\Omega'' \setminus \bar{\Omega}') 
         \leq C(\eta) \epsilon.
         \]
         Putting these inequalities together, it holds
         \begin{align*}
    P(\Epi(\phi+u),\Omega'' \times \bb{R}) \geq P(\Epi(f),\Omega'' \times \bb{R})  
         -\epsilon
         \\
         \geq P(\Epi(f \eta + (1-\eta)u),\Omega'' \times \bb{R})-(C(\eta)+1)\epsilon
         \\
         \geq 
         P(\Epi(u),\Omega'' \times \bb{R}) -(C(\eta)+1))\epsilon,
         \end{align*}
         concluding the proof.
    \end{proof}
\end{proposition}

Given $f :\X \times \bb{R} \to \bb{R}$ and $(x,t) \in \X \times \bb{R}$ we denote by $f^t$ and $f^x$ respectively the restriction of $f$ to $\X \times \{t\}$ and to $\{x\} \times \bb{R}$.

\begin{thm}
    Let $(\X,\sd,\m)$ be an $\RCD(K,N)$ space and let $\Omega \subset \X$ be open. Let $u \in \W^{1,1}(\Omega)$ be a solution of the minimal surface equation. Then $\Epi(u)$ is a perimeter minimizing set in $\Omega \times \bb{R}$.
    \begin{proof}
    Consider a set $E \subset \Omega \times \bb{R}$ such that $\Epi'(u) \Delta E \subset \subset \Omega \times \bb{R}$. Modulo translating vertically, we may suppose that there exists $c>1$ such that $u$ takes values in $(1,c-1)$ and $\Epi'(u) \Delta E \subset \subset \Omega \times (1,c-1)$.
    We then define $w(E):\Omega \to \bb{R}$ by
    \[
    w(E)(x):=\int_0^c 1_E(x,s) \, ds
    \]
    and we claim that
    \begin{equation} \label{E|w(E)}
    w(E) \in \BV(\Omega)
    \quad 
    \text{and}
    \quad
    P(\Epi'(w(E)),\Omega \times \bb{R}) \leq P(E,\Omega \times \bb{R}).
    \end{equation}
    If the claim holds, noting that by construction $\{w(E) \neq u\} \subset \subset \Omega$, by Proposition \ref{P22}, we get
    \[
    P(\Epi'(u),\Omega \times \bb{R}) \leq P(\Epi'(w(E)),\Omega \times \bb{R})\leq P(E,\Omega \times \bb{R}).
    \]
    To prove the claim, we consider a sequence $f_n \in \Lip(\Omega \times (0,c))$ converging in $\sL^1(\Omega \times (0,c))$ to $1_E$ and such that
    \[
    \lim_{n \to + \infty} |D f_n|(\Omega \times (0,c))=P(E,\Omega \times \bb{R}).
    \]
    Modulo truncating, we can assume that $f_n \equiv 1$ on $\Omega \times \{0\}$ and $f_n \equiv 0$ on $\Omega \times \{c\}$.
    We then define $w(f_n): \Omega \to \bb{R}$ as
    \[
    w(f_n)(x):=\int_0^c f_n(x,s) \, ds.
    \]
    These functions are Lipschitz since
    \[
    \frac{|w(f_n)(x)-w(f_n)(y)|}{\sd(x,y)} \leq \int_0^c \frac{|f_n(x,s)-f_n(y,s)|}{\sd(x,y)} \, ds \leq c \ssf{L}(f_n).
    \]
    By reverse Fatou Lemma and the fact that each $f_n$ is Lipschitz, for $\m$-a.e. $x \in \X$ it holds
    \begin{equation} \label{E30}
    |\nabla w(f_n)|(x)=\lip (w(f_n))(x) \leq \int_0^c |\nabla f_n^s|(x) \, ds.
    \end{equation}
    Moreover, by the tensorization property of the Cheeger energy (see \cite{Tens}), it holds
    \begin{align*}
    |D f_n|(\Omega \times (0,c))
    & =
    \int_{\Omega \times (0,c)} |\nabla f_n |(x,s) \, d\m \, ds =
    \int_{\Omega \times (0,c)} \sqrt{|\nabla f_n^s|^2(x)+|\nabla f_n^x|^2(s)} \, d\m \, ds \\
    & \geq 
    \sup_{\substack{(a,b) \in C(\Omega) \times C(\Omega) \\ a^2+b^2 \leq 1 \\ a,b \geq 0}}
    \int_{\Omega \times (0,c)} a(x) |\nabla f_n^s|(x)+b(x)|\nabla f_n^x|(s) \, d\m \, ds.
    \end{align*}
    Combining with \eqref{E30} and the fact that $f_n=1$ on $\Omega \times \{0\}$ and $f_n=0$ on $\Omega \times \{c\}$, we deduce
    \begin{align*}
        |D f_n|(\Omega \times (0,c)) \geq 
    \sup_{\substack{(a,b) \in C(\Omega) \times C(\Omega) \\ a^2+b^2 \leq 1 \\ a,b \geq 0}}
    \int_\Omega a(x) |\nabla w(f_n)|(x)+b(x) \, d\m=\int_\Omega \sqrt{1+|\nabla w(f_n)|^2} \, d\m.
    \end{align*}
    Moreover, $w(f_n) \to w(E)$ in $\sL^1(\Omega)$ as $n \to + \infty$ since
    \[
    \int_\Omega |w(f_n) - w(E)| \, d\m \leq \int_{\Omega \times (0,c)}|f_n(x,s)-1_E(x,s)| \, d\m \, ds.
    \]
    These facts imply that $w(E) \in \BV(\Omega)$ as claimed. Concerning the perimeter of $w(E)$, we get
    \begin{align*}
    P(\Epi'(w(E)),\Omega \times \bb{R}) & \leq \liminf  P(\Epi'(w(f_n)),\Omega \times \bb{R}) \\
    & \leq \liminf |D f_n|(\Omega \times (0,c))=P(E,\Omega \times \bb{R}),
    \end{align*}
    proving the claim \eqref{E|w(E)}.
    \end{proof}
\end{thm}

\section{Harmonicity of the height function on a minimal graph} \label{S2}

In this section, we prove Theorems \ref{CT5} and \ref{CT6} from the introduction. In Section \ref{S8}, we study the blow-up properties of solutions of the minimal surface equation on $\RCD(K,N)$ spaces. In Section \ref{S9}, we equip minimal graphs on $\RCD(K,N)$ spaces with a metric measure structure, and we prove Theorems \ref{CT5} and \ref{CT6} from the introduction.

\subsection{Blow-ups of functions} \label{S8}

In \cite[Theorem 10.2]{Cheeger}, it is shown that a Lipschitz function $f$ on a PI space $(\X,\sd,\m)$ is infinitesimally generalized linear for $\m$-almost every $x$. This means that, when fixing a tangent cone to $\X$ at $x$, the corresponding blow-up of $f$ at $x$ is a function $f^\infty$ which is linear in a generalized sense. Moreover, the local Lipschitz constant of $f$ at $x$ and the Lipschitz constant of $f^\infty$ coincide.

In this section, we consider a variation of this property. Given an $\RCD(K,N)$ space $(\X,\sd,\m)$ and a solution of the minimal surface equation $u \in \W^{1,1}_{loc}(\Omega)$, we would like to use an analogue of \cite[Theorem 10.2]{Cheeger} for $u$. Since $u$ is not Lipschitz, we cannot apply Cheeger's result directly. 
Theorem \ref{T16} shows that $u$ is nevertheless infinitesimally generalized linear at $\m$-almost every point of its domain.

The next two definitions are needed to exploit the machinery from \cite{Cheeger} in our context.

\begin{definition} \label{D3}
    Let $(\X,\sd,\m)$ be an $\RCD(K,N)$ space, let $x \in \X$, and let $\Tan_x(\X)=\{(\bb{R}^k,\sd_e)\}$ (i.e. $x$ is a regular point of $\X$). Given a sequence $\{\epsilon_n\}_{n \in \bb{N}} \subset (0,+\infty)$ such that $n \epsilon_n \to 0$, we say that $(\epsilon_n,\psi_n)$ is a blow-up of $\X$ at $x$ if there exists a sequence $\delta_n$ with $n \delta_n \to 0$ such that each
    $\psi_n:\bar{B}_n(0) \subset \bb{R}^k \to (\bar{B}_{n\epsilon_n}(x), \epsilon_n^{-1}\sd)$ is a $\delta_n$-$GH$ map.
\end{definition}

We stress that for any point $x \in \X$ with $\Tan_x(\X)=\{(\bb{R}^k,\sd_e)\}$ and any sequence $\epsilon_n \to 0$, there exists a (non relabeled) subsequence such that there exists a blow-up $(\epsilon_n,\psi_n)$ of $\X$ at $x$.

Given a bounded function $f:\X \to \bb{R}$, we denote by $\|f\|_{\infty,\X}$ the supremum of $|f|$. 

\begin{definition} \label{D8}
    Let $(\X,\sd,\m)$ be an $\RCD(K,N)$ space, let $x \in \X$ be a regular point, and let $f:\X \to \bb{R}$ be a function. Given a sequence $\{\epsilon_n\}_{n \in \bb{N}}$ with $n \epsilon_n \to 0$, we say that a triple $(\epsilon_n,\psi_n,f^{\infty})$ is a blow-up of $f$ at $x$ if $(\epsilon_n,\psi_n)$ is a blow-up of $\X$ at $x$ and
    $f^{\infty}: \bb{R}^k \to \bb{R}$ is a linear function such that
        \[
        \|\epsilon_n^{-1}(f \circ \psi_n - f(x))-f^{\infty}\|_{\infty,\bar{B}^{\bb{R}^k}_1(0)} \to 0.
        \]
        If, for every blow-up $(\epsilon_n,\psi_n)$ of $\X$ at $x$, the function $f$ admits a blow-up at $x$ (depending on $(\epsilon_n,\psi_n)$), we say that $f$ is infinitesimally generalized linear at $x$.
\end{definition}

We remark that the definition given in \cite[Definition 10.1]{Cheeger} of a point $x \in \X$ where a Lipschitz function $f \in \Lip(\X)$ is infinitesimally generalized linear requires in addition that that $x$ is a Lebesgue point for $\lip(f)$.

\begin{remark} \label{R2}
    Let $(\X,\sd,\m)$ be an $\RCD(K,N)$ space, let $x \in \X$ be a regular point, and let $f:\X \to \bb{R}$ be an infinitesimally generalized linear function at $x$. Let $(\epsilon_n,\psi_n)$ and $(\epsilon_n,\psi'_n)$ be blow-ups of $\X$ at $x$ such that $\psi_n \neq \psi'_n$ on at most a finite number of points. Let $(\epsilon_n,\psi_n,f^\infty)$ and $(\epsilon_n,\psi'_n,{f'}^\infty)$ be the corresponding blow-ups of $f$. Then $f^\infty$ and ${f'}^\infty$ are linear functions coinciding outside of a countable set of $\bb{R}^k$, so that $f^\infty = {f'}^\infty$.
\end{remark}


The next theorem follows from \cite[Theorem $10.2$]{Cheeger} and Theorem \ref{constancy}, and it is the key tool of this section. 

\begin{thm} \label{T3}
    Let $(\X,\sd,\m)$ be an $\RCD(K,N)$ space, let $\Omega \subset \X$ be an open set and let $f \in \Lip(\Omega)$. Then for $\m$-a.e. $x \in \Omega$, for every blow-up $(\epsilon_n,\psi_n)$ of $\X$ at $x$ there exists a (non relabeled) subsequence and a linear function $f^\infty:\bb{R}^k \to \bb{R}$ such that $(\epsilon_n,\psi_n,f^{\infty})$ is a blow-up of $f$ at $x$ and
    $| \nabla f|(x)=\lip(f)(x)=\lip(f^{\infty})(0)$.
\end{thm}

The next lemma and the subsequent proposition are the key tools to prove Theorem \ref{T16}. 
Lemma \ref{L20} shows the following. Let $(\X,\sd,\m)$ be an $\RCD(K,N)$ space and let $x \in \X$ be a regular point of $\X$. Let $E$ be a perimeter minimizing set admitting tangent balls at $x \in \partial E$. Let $(\epsilon_n,\psi_n)$ be a blow-up of $\X$ in $x$. There are pairs of balls of arbitrarily large radius in $\bb{R}^k$ (whose centers lie on a line through the origin, and whose radius is approximately the distance of their centers from the origin) which are sent by $\psi_n$, for $n$ large enough, respectively in $E$ and $\X \setminus E$.

\begin{lemma} \label{L20}
    Let $(\X,\sd,\m)$ be an $\RCD(K,N)$ space and let $x \in \X$ be a regular point of $\X$. Let $E$ be a perimeter minimizing set admitting tangent balls at $x \in \partial E$. Let $(\epsilon_n,\psi_n)$ be a blow-up of $\X$ in $x$. There exist $\{a_n\}_{n \in \bb{N}} \subset \bb{R}^k$, $\{c_n\}_{n \in \bb{N}} \subset \bb{R}_+$, with $|a_n| \to + \infty$ and $c_n \to 0$, such that
    \[
    \psi_n(B^{\bb{R}^k}_{|a_n|-c_n}(a_n)) \subset E \setminus \partial E, \quad \psi_n(B^{\bb{R}^k}_{|a_n|-c_n}(-a_n)) \subset \X \setminus (E \cup \partial E).
    \]
    \begin{proof}
        Let $\gamma$ be a geodesic connecting the centers of the tangent balls to $E$ at $x$ which passes through $x$. We remark that such a geodesic exists, since  concatenating geodesics from
centers of tangent balls to the corresponding tangent point still yields a geodesic.
        
        Let $y_1,y_2 \in \partial B_{n\epsilon_n}(x)$ be the points of intersection of $\gamma$ with $\partial B_{n\epsilon_n}(x)$. 

        Since the maps $\psi_n$ are $\delta_n$-GH maps (for some $\delta_n>0$ such that $n \delta_n \to 0$), there exist points $b^n_1,b^n_2 \in B_n^{\bb{R}^k}(0)$ such that $\epsilon_n^{-1}\sd(\psi_n(b^n_1),y_1) \leq \delta_n$ and $\epsilon_n^{-1}\sd(\psi_n(b^n_2),y_2) \leq \delta_n$. Hence,
        $
        |b^n_1-b^n_2| \geq 2n-3\delta_n,
        $
        so that
        \begin{equation} \label{E58}
        |b^n_1+b^n_2|^2 = 2|b^n_1|^2+2|b^n_2|^2-|b_1^n-b^n_2|^2 \leq 12n \delta_n+9\delta_n^2. 
        \end{equation}
        Setting $a_n:=b^n_1$, it trivially holds
        $
        \epsilon_n^{-1}
        \sd(\psi_n(a_n),y_1) \leq \delta_n
        $. At the same time, using \eqref{E58} in the last inequality, it holds
        \[
        \epsilon_n^{-1}
        \sd(\psi_n(-a_n),y_2) \leq \epsilon_n^{-1}
        \sd(\psi_n(b_2^n),y_2)
        + |b^n_2+a_n|+\delta_n \leq 12n \delta_n+9\delta_n^2+2\delta_n.
        \]
        We set $c_n:=100(12n \delta_n+9\delta_n^2+2\delta_n)$, so that 
        \[
        \epsilon_n^{-1}
        \sd(\psi_n(-a_n),y_2) \leq c_n/100, \quad \epsilon_n^{-1}
        \sd(\psi_n(a_n),y_1) \leq c_n/100.
        \]
        Using again that $\psi_n$ almost preserves distances, we deduce
        \[
        \psi_n(B_{n-c_n}(a_1)) \subset B^\X_{n \epsilon_n}(y_1) \, \quad \psi_n(B_{n-c_n}(-a_1)) \subset B^\X_{n \epsilon_n}(y_2).
        \]
        Since $n \epsilon_n \to 0$, for $n$ large enough, the balls $B^\X_{n \epsilon_n}(y_2), B^\X_{n \epsilon_n}(y_1) $ are contained in the tangent balls to $E$ at $x$. Since $n \delta_n \to 0$, it holds $c_n \to 0$, concluding the proof.
    \end{proof}
\end{lemma}

The next proposition shows that the graph of a solution of the minimal surface equation on an $\RCD(K,N)$ space is well approximated at small scales by a hyperplane at points where the epigraph admits tangent balls.

\begin{proposition} \label{P34}
    Let $(\X,\sd,\m)$ be an $\RCD(K,N)$ space, let $\Omega \subset \X$ be an open set, and let $u \in \W^{1,1}(\Omega)$ be a solution of the minimal surface equation. 
    Let $x \in \Omega$ be a regular point of $\X$ such that $\Epi(u)$ has tangent balls at $(x,u(x))$. Let $(\epsilon_n,\psi_n)$ be a blow-up of $\X$ at $x$. Modulo passing to a subsequence, there exists an hyperplane $H \subset \bb{R}^k \times \bb{R}$ satisfying the following. Let $v:\Omega \to \bb{R}$ be a function with $\Graph(v) \subset \partial \Epi(u)$ and $v(x)=u(x)$.
    \begin{enumerate}
        \item \label{iitem1} If $y_n \in \Graph(\epsilon_n^{-1} v \circ \psi_n-v(x)) \subset \bb{R}^k \times \bb{R}$ converge to $y$, then $y \in H$.
    \item \label{iitem2} If $H$ is the graph of a linear function $f:\bb{R}^k \to \bb{R}$, then
    \[
    \|\epsilon_n^{-1}(v \circ \psi_n-v(x))-f\|_{\infty,\bar{B}_1^{\bb{R}^k}(0)} \to 0.
    \]
    \end{enumerate}
    \begin{proof}
        Modulo translating, we may assume that $u(x)=0$.
        Consider the maps $\tilde{\psi}_n:=(\psi_n,\epsilon_n \ssf{Id}): B_n^{\bb{R}^k \times \bb{R}}(0) \to (\X \times \bb{R},\epsilon_n^{-1}\sd_\times,(x,0))$ given by $\tilde{\psi}_n(y,t):=(\psi_n(y),\epsilon_nt)$. By construction, $(\epsilon_n,\tilde{\psi}_n)$ is a blow-up of $\X \times \bb{R}$ in $(x,u(x))$.

        By Lemma \ref{L20}, there exist
        $\{a_n\}_{n \in \bb{N}} \subset \bb{R}^{k+1}$, $\{c_n\}_{n \in \bb{N}} \subset \bb{R}_+$, with $|a_n| \to + \infty$ and $c_n \to 0$, such that
    \[
    \tilde{\psi}_n
    (B^{\bb{R}^{k+1}}_{|a_n|-c_n}(a_n)) \subset \Epi(u) \setminus \partial \Epi(u), \quad \tilde{\psi}_n(B^{\bb{R}^{k+1}}_{|a_n|-c_n}(-a_n)) \subset \X \times \bb{R} \setminus (\Epi(u) \cup \partial \Epi(u)).
    \]
    Passing to a subsequence, the sets 
    \[
    \bb{R}^{k+1} \setminus (B_{|a_n|-c_n}(-a_n) \cup B_{|a_n|-c_n}(a_n))
    \]
    converge in Hausdorff sense in any compact set of $\bb{R}^{k+1}$ to an hyperplane $H$.
    
    Let $v:\Omega \to \bb{R}$ be as in the statement. 
    To prove item \ref{iitem1}, by our previous remark, it is sufficient to show that 
    \begin{equation} \label{E60}
    \Graph(\epsilon_n^{-1} v \circ \psi_n) \cap (B_{|a_n|-c_n}(-a_n) \cup B_{|a_n|-c_n}(a_n))=\emptyset.
    \end{equation}
    If \eqref{E60} fails, then, without loss of generality, there exists $y \in \bb{R}^k$ such that $(y,\epsilon_n^{-1}v(\psi_n(y))) \in B_{|a_n|-c_n}(a_n)$. Hence, 
    \[
    (\psi_n(y),v(\psi_n(y)))=\tilde{\psi_n}(y,\epsilon_n^{-1}v(\psi_n(y))) \in \tilde{\psi}_n
    (B^{\bb{R}^{k+1}}_{|a_n|-c_n}(a_n)),
    \]
    so that $(\psi_n(y),v(\psi_n(y))) \in \Epi(u) \setminus \partial \Epi(u)$, a contradiction.

    To prove item \ref{iitem2}, it is sufficient to show that 
    \begin{align} \label{E61}
    & \Graph(\epsilon_n^{-1} v \circ \psi_n) \cap \{(a,t) \in \bb{R}^k \times \bb{R}: \exists \, s<t: (a,s) \in B_{|a_n|-c_n}(a_n) \} =\emptyset \\
    & \nonumber 
    \Graph(\epsilon_n^{-1} v \circ \psi_n) \cap \{(a,t) \in \bb{R}^k \times \bb{R}: \exists \, s>t: (a,s) \in B_{|a_n|-c_n}(-a_n) \}
    =\emptyset.
    \end{align}
    Assume that \eqref{E61} fails, the other case being analogous.
    Then, there exists $y \in \bb{R}^k$ and $s<\epsilon_n^{-1}v(\psi_n(y))
    $
    such that $(y,s) \in B_{|a_n|-c_n}(a_n)$. As for item \ref{iitem1}, we deduce $(\psi_n(y),\epsilon_ns) \in \Epi(u) \setminus \partial \Epi(u)$, so that by definition of epigraph it also follows $(\psi_n(y),v(\psi_n(y)) \in \Epi(u) \setminus \partial \Epi(u)$, a contradiction.
    \end{proof}
\end{proposition}

The next theorem, which is needed for Theorem \ref{T16}, can be found in \cite[Theorem $6.12$]{EvansGari}. The proof in the $\RCD(K,N)$ setting is analogous.
\begin{thm} \label{T14}
    Let $(\X,\sd,\m)$ be an $\RCD(K,N)$ space and let $f \in \BV(\X)$. Then, for every $\epsilon>0$, there exists a Lipschitz function $f_\epsilon$ such that $\m(\{f \neq f_\epsilon \}) \leq \epsilon$.
\end{thm}

We now prove the main theorem of the section, showing that 
a solution $u$ of the minimal surface equation is infinitesimally generalized linear at $\m$-almost every point of its domain.
The idea is the following. 
By the previous theorem, given a Lipschitz function $f$, it is enough to show the statement at points where $\{f=u\}$ has density $1$. So we restrict ourselves to these points. By Cheeger's Theorem \ref{T3}, $f$ will have a blow-up given by a linear function $f^\infty:\bb{R}^k \to \bb{R}$ at $x$. By the previous proposition, modulo staying out of a negligible set, the graph of $u$ at $(x,u(x))$ is well approximated at small scales by a Euclidean hyperplane $H \subset \bb{R}^{k+1}$. Since we were at a density point of $\{f=u\}$, this forces $H=\Graph(f^\infty)$, giving the desired conclusion.

\begin{thm} \label{T16}
    Let $(\X,\sd,\m)$ be an $\RCD(K,N)$ space, let $\Omega \subset \X$ be an open set, and let $u \in \W^{1,1}(\Omega)$ be a solution of the minimal surface equation. For $\m$-almost every regular point $x \in \Omega$, the following happens. Let $(\epsilon_n,\psi_n)$ be a blow-up of $\X$ at $x$. Modulo passing to a subsequence, there exists a linear function $u^\infty:\bb{R}^k \to \bb{R}$ such that any $v:\Omega \to \bb{R}$ with $\Graph(v) \subset \partial \Epi(u)$ has blow-up $(\epsilon_n,\psi_n,u^\infty)$ at $x$.
    \begin{proof}
        Let $\epsilon>0$ and let $f \in \Lip(\Omega)$ be such that $\m(\{f \neq u\}) \leq \epsilon$. It is sufficient to show the statement at a regular point $x \in \X$ where $f$ is infinitesimally generalized linear and $\{f=u\}$ has density $1$. Modulo translating, we may assume $u(x)=f(x)=0$.

        Let $(\epsilon_n,\psi_n)$ be a blow-up of $\X$ at $x$.
        Let $r>0$ and let $x_1,...,x_k$ be points in $\bb{R}^k$ such that each ball $\bar{B}_r(x_i)$ is contained in $B^{\bb{R}^k}_1(0)$, and such that any $k$-tuple $(p_1,\cdots,p_k)$, such that $p_i \in \bar{B}_r(x_i)$, is made of linearly independent points.

        Since $\{f=u\}$ has density $1$ in $x$, for $n$ large enough, $B_{r/2}(\psi_n (x_i)) \cap \{f=u\}$ is non-empty for every $i=1, \cdots k$.
        Hence, we can consider a blow-up $(\epsilon_n,\psi'_n)$ of $\X$ at $x$, such that $\psi'_n \neq \psi_n$ at a finite number of points and, for every $i \in \{1,\cdots,k\}$, there exists $p^n_i \in B_r(x_i)$ with $\psi'_n (p^n_i) \in \{f=u\}$. The $k$-tuples $(p^n_1, \cdots,p^n_k)$ converge to a $k$-tuple $(p_1,\cdots ,p_k)$ of $k$ linearly independent points.
        

        Since $f$ is Lipschitz, by Theorem \ref{T3}, there exists a linear function $f^\infty:\bb{R}^k \to \bb{R}$ such that
        \[
        \|\epsilon_n^{-1}f \circ \psi'_n-f^\infty\|_{\infty,B_1^{\bb{R}^k}(0)} \to 0.
        \]
        Modulo passing to a subsequence, let $H' \subset \bb{R}^{k+1}$ be the hyperplane relative to $(\epsilon_n,\psi'_n)$ given by Proposition \ref{P34}. Since by construction the points $\{(p_i,f^\infty(p_i))\}_{i=1}^k$ belong to $H'$, it holds $H'=\Graph(f^\infty)$.

        Let now $H \subset \bb{R}^{k+1}$ be the hyperplane relative to $(\epsilon_n,\psi_n)$ given by Proposition \ref{P34}. Out of the countable set where $\psi_n \neq \psi'_n$ for some $n$, the functions $\epsilon_n^{-1}u \circ \psi_n$ coincide with $\epsilon^{-1}_nu \circ \psi'_n$ and converge pointwise to $f^\infty$. Hence $H=H'=\Graph(f^\infty)$, so that the statement follows by item \ref{iitem2} in Proposition \ref{P34}.

    \end{proof}
\end{thm}

\subsection{Minimal graphs as metric measure spaces} \label{S9}

In this section, we prove Theorems \ref{CT5} and \ref{CT6} from the introduction. 

Let $(\X,\sd,\m)$ be an $\RCD(K,N)$ space and let $u \in \W^{1,1}_{loc}(\Omega)$ be a solution of the minimal surface equation. 
We would like to say that Lipschitz functions on $\X \times \bb{R}$, when restricted to the graph of $u$, are infinitesimally generalized linear at almost every point. Since we do not know whether the graph of $u$ is a PI space, we cannot apply Cheeger's result \cite[Theorem 10.2]{Cheeger} directly. This difficulty is overcome with Theorem \ref{T20}, which in turn relies on the results from the previous section.

We first define the metric measure space structure on the graph of a solution of the minimal surface equation.

\begin{definition} \label{D5}
    Let $(\X,\sd,\m)$ be an $\RCD(K,N)$ space, let $\Omega \subset \X$ be open, and let $u \in \W^{1,1}_{loc}(\Omega)$ be a solution of the minimal surface equation. Let $\G(u) \subset \Omega \times \bb{R}$ be the boundary of $\Epi(u)$. We define the metric measure space $(\G(u),\sd_g,\m_g)$, where $\sd_g$ is the product distance of $\Omega \times \bb{R}$ restricted to $\G(u) \times \G(u)$ and $\m_g$ is the perimeter measure induced by $\Epi(u)$ on the Borel subsets of $\G(u)$. 
\end{definition}

We denote by $\lip_g$ the local Lipschitz constant of a function $f:\G(u) \to \bb{R}$ and by $i_g:\Omega \to \G(u)$ the function given by \begin{equation} \label{E71}
    i_g(x):=(x,u(x)).
\end{equation}

We now prove a few short technical lemmas describing the behavior of Lipschitz functions on $\G(u)$.

\begin{lemma} \label{L15}
    Let $(\X,\sd,\m)$ be an $\RCD(K,N)$ space, let $\Omega \subset \X$ be an open set, let $u \in \W_{loc}^{1,1}(\Omega) \cap \sL^\infty_{loc}(\Omega)$, and let $f \in \Lip(\X \times \bb{R})$. Let $g:\Omega \to \bb{R}$ be defined by $g(x):=f(x,u(x))$. Then, $g \in \W_{loc}^{1,1}(\Omega)$.
    \begin{proof}
        Modulo multiplying by a cut off function, we may suppose that $u \in \W^{1,1}(\X) \cap \sL^\infty(\X)$ with compact support. Let $u_n \in \Lip_{loc}(\X) \cap \sL^1(\X)$ be such that $u_n \to u$ in $\sL^1(\X)$, with uniformly bounded norms of the gradients $\|\nabla u_n\|_{\sL^1(\X)}$, and such that $|\nabla u_n| \, d\m \to |\nabla u| \, d\m$ weakly in duality with continuous compactly supported functions. 

        Consider the function $g_n \in \Lip_{loc}(\X)$ defined by $g_n(x):=f(x,u_n(x))$. It is easy to check that $\lip(g_n) \leq \ssf{L}(f)(1+|\nabla u_n|)$. Since $g_n \to g$ in $\sL^1_{loc}(\X)$, it holds $g \in \BV_{loc}(\X)$ and $|D g| \leq \ssf{L}(f)(1+|\nabla u|)$. In particular, $g \in \W^{1,1}_{loc}(\X)$.
    \end{proof}
\end{lemma}

\begin{lemma} 
    Let $(\X,\sd,\m)$ be an $\RCD(K,N)$ space, let $\Omega \subset \X$ be an open set, and let $u \in \W^{1,1}(\Omega)$. For $\m$-a.e. regular point of $x \in \X$ where $u$ is infinitesimally generalized linear, if $(\epsilon_n,\psi_n,u^\infty)$ is a blow-up of $u$ at $x$, then $|\nabla u|(x)=\lip(u^\infty)$.
    \begin{proof}
        Let $\epsilon>0$ and let $f \in \Lip(\Omega)$ be a function with $\m(\{f \neq u\}) \leq \epsilon$. Let $x$ be a point where $u$ and $f$ are infinitesimally generalized linear, with $|\nabla u|(x)=|\nabla f|(x)$, and assume that $\{f = u\}$ has density $1$ at $x$.

        Let $(\epsilon_n,\psi_n,u^\infty)$ and $(\epsilon_n,\psi_n,f^\infty)$ be blow-ups of $u$ and $f$ respectively. 
        To conclude, by Theorem \ref{T3}, it is enough to show that $u^\infty=f^\infty$. Since we are at a point where both $f$ and $u$ are infinitesimally generalized linear, if we consider a blow-up $(\epsilon_n,\psi'_n)$ of $\X$ at $x$, such that $\psi_n \neq \psi'_n$ at a finite number of points, then $(\epsilon_n,\psi'_n,u^\infty)$ and $(\epsilon_n,\psi'_n,f^\infty)$ are still blow-ups of $u$ and $f$.
        Moreover, since $f^\infty$ and $u^\infty$ are linear, it is sufficient to show that they coincide in $k$ linearly independent points of $\bar{B}^{\bb{R}^k}_1(0)$.

        Let $r>0$ and let $x_1,...,x_k$ be points in $\bb{R}^k$ such that each ball $\bar{B}_r(x_i)$ is contained in $B^{\bb{R}^k}_1(0)$, and such that any $k$-tuple $(p_1,\cdots,p_k)$, with that $p_i \in \bar{B}_r(x_i)$, is made of linearly independent points.

        Since $\{f=u\}$ has density $1$ in $x$, for $n$ large enough, $B_{r/2}(\psi_n(x_i)) \cap \{f=u\}$ is non-empty for every $i=1, \cdots k$.
        Hence, we can consider a blow-up $(\epsilon_n,\psi'_n)$ of $\X$ at $x$, such that $\psi'_n \neq \psi_n$ at a finite number of points and, for every $i \in \{1,\cdots,k\}$, there exists $p^n_i \in B_r(x_i)$ with $\psi'_n (p^n_i) \in \{f=u\}$. The $k$-tuples $(p^n_1, \cdots,p^n_k)$ converge to a $k$-tuple $(p_1,\cdots ,p_k)$ of $k$ linearly independent points. Hence, as $n \to + \infty$, it holds
        \begin{align*}
        & |f^\infty (p_i)-  u^\infty(p_i)| \leq |f^\infty (p_i)-f^\infty(p^n_i)|+|u^\infty(p_i)-u^\infty(p^n_i)|  \\
        & +
        |\epsilon_n^{-1}(f(\psi'_n (p^n_i))-f(x))-f^\infty(p^n_i)|+|\epsilon_n^{-1}(u(\psi'_n (p^n_i))-u(x))-u^\infty(p^n_i)|
        \to 0.
        \end{align*}
        This shows that $f^\infty (p_i)= u^\infty(p_i)$ for every $i \in \{1,\cdots ,k\}$, concluding the proof.
    \end{proof}
\end{lemma}

\begin{corollary} \label{CL19}
    Let $(\X,\sd,\m)$ be an $\RCD(K,N)$ space, let $\Omega \subset \X$ be an open set, and let $u,v \in \W^{1,1}(\Omega)$. For $\m$-a.e. regular point of $x \in \X$ where both $u$ and $v$ are infinitesimally generalized linear, the following happens. If $(\epsilon_n,\psi_n,u^\infty)$ and $(\epsilon_n,\psi_n,v^\infty)$
    are blow-ups of $u$ and $v$ at $x$, then $\nabla u \cdot \nabla v=\nabla u^\infty \cdot \nabla v^\infty$.
\end{corollary}

Let $(\X,\sd,\m)$ be an $\RCD(K,N)$ space, let $\Omega \subset \X$ be an open set, let $u \in \W^{1,1}(\Omega)$ be a solution of the minimal surface equation.
The next lemma shows that, for almost every point of $\G(u)$, the local Lipschitz constant of a function on $\G(u)$ coincides with the local Lipschitz constant of its restriction to $\Graph(u)$. We recall that $u$ is continuous $\m$-almost everywhere by Proposition \ref{P14}.

\begin{lemma} \label{L17}
    Let $(\X,\sd,\m)$ be an $\RCD(K,N)$ space, let $\Omega \subset \X$ be an open set, let $u \in \W^{1,1}(\Omega)$ be a solution of the minimal surface equation. Let $\theta \in \Lip(\G(u))$ and let $x \in \Omega$ be a point where $u$ is infinitesimally generalized linear and continous. Then, there exist points $x_n \to x$ such that 
    \[
        \lip_g(\theta)(i_g(x))=\lim_{n \to + \infty} \frac{|\theta(i_g(x_n))-\theta(i_g(x))|}{\sd_g(i_g(x),i_g(x_n))}.
        \]
        \begin{proof}
            Let $\bar{x}_n:=(x_n,t_n) \in \G(u)$ be points such that
            \[
        \lip_g(\theta)(i_g(x))=\lim_{n \to +\infty} \frac{|\theta(\bar{x}_n)-\theta(i_g(x))|}{\sd_g(i_g(x),\bar{x}_n)}.
        \]
        Since $u$ is continuous at $x$, it holds $\{x\} \times \bb{R} \cap \G(u)=\{(x,u(x))\}$. In particular, $x_n \neq x$.
        It is sufficient to show that
        \[
        \lip_g(\theta)(i_g(x)) \leq \limsup_{n \to + \infty} \frac{|\theta(i_g(x_n))-\theta(i_g(x))|}{\sd_g(i_g(x),i_g(x_n))}.
        \]
        We first show
        \begin{equation} \label{E48}
        \frac{|t_n-u(x_n)|}{\sd(x,x_n)} \to 0.
        \end{equation}
        This then implies
        \begin{equation} \label{E49}
            \lim_{n \to + \infty}
        \frac{\sd_g(i_g(x),\bar{x}_n)}{
        \sd_g(i_g(x),i_g(x_n))}=1.
        \end{equation}
        Let $\epsilon_n:=\sd(x,x_n)$ and let $(\epsilon_n,\psi_n)$ be a blow-up of $\X$ at $x$ such that $x_n \in \ssf{Im}(\psi_n)$ with $\psi_n(x'_n)=x_n$. By Theorem \ref{T16}, there exists a linear function $u^\infty: \bb{R}^k \to \bb{R}$ such that $(\epsilon_n,\psi_n,u^\infty)$ is a blow-up of any function $v:\Omega \to \bb{R}$ with $\Graph(v) \subset \G(u)$. Hence, setting $v=u$ outside of $\{x_n\}$ and $v(x_n):=t_n$, it holds
        \[
        \frac{|t_n-u(x_n)|}{\sd(x,x_n)} \leq |\epsilon_n^{-1}(v(\psi_n(x'_n))-v(x))-u^\infty(x'_n)|+|\epsilon_n^{-1}(u(\psi_n(x'_n))-u(x))-u^\infty(x'_n)|.
        \]
        Since the r.h.s. converges to zero by the previous remarks, $\frac{|t_n-u(x_n)|}{\sd(x,x_n)} \to 0$ as well, proving \eqref{E48}.

        By \eqref{E49}, to conclude, we need to show that
\begin{equation} \label{E56}
\lim_{n \to + \infty} \frac{|\theta(\bar{x}_n)-\theta(i_g(x))|}{\sd_g(i_g(x_n),i_g(x))} \leq 
\limsup_{n \to + \infty} \frac{|\theta(i_g(x_n))-\theta(i_g(x))|}{\sd_g(i_g(x_n),i_g(x))}.
\end{equation}
Since $\theta$ is Lipschitz and $\sd_g(i_g(x_n),i_g(x)) \geq \sd(x_n,x)$, it holds
\begin{align*}
\frac{|\theta(\bar{x}_n)-\theta(i_g(x))|}{\sd_g(i_g(x_n),i_g(x))}
\leq \frac{|\theta(i_g(x_n))-\theta(i_g(x))|}{\sd_g(i_g(x_n),i_g(x))}
+ \frac{\ssf{L}(\theta)|t_n-u(x_n)|}{\sd(x_n,x)}.
\end{align*}
Combining with \eqref{E48}, we deduce \eqref{E56}.
        \end{proof}
\end{lemma}

The next lemma shows a locality property for the local lipschitz constant of a function on a minimal graph.

\begin{lemma} \label{L16}
    Let $(\X,\sd,\m)$ be an $\RCD(K,N)$ space, let $\Omega \subset \X$ be an open set, let $u \in \W^{1,1}(\Omega)$ be a solution of the minimal surface equation. Let $x \in \Omega$ be a point where $u$ is infinitesimally generalized linear and continuous. Let $\theta_1,\theta_2 \in \Lip(\G(u))$ be such that $\theta_1(i_g(x))=\theta_2(i_g(x))$ and $\{\theta_1 \circ i_g=\theta_2 \circ i_g\}$ has density $1$ in $x$. Then, $\lip_g(\theta_1)(x)=\lip_g(\theta_2)(x)$.
    \begin{proof}
        It is sufficient to show that $\lip_g(\theta_1-\theta_2)=0$. Let $f:=\theta_1-\theta_2$. Let $x_n \to x$ be the sequence given by the previous lemma such that
        \[
        \lip_g(f)(x)=\lim_{n \to + \infty} \frac{|f(i_g(x))-f(i_g(x_n))|}{\sd_g(i_g(x_n),i_g(x))}.
        \]
        Since $\{f \circ i_g=0\}$ has density $1$ in $x$, for every $\epsilon>0$ and $N \in \bb{N}$, there exist $n \geq N$ and $y_n \in \{f \circ i_g=0\}$ such that $\sd(y_n,x_n) \leq \epsilon \sd(x,x_n)$. Since $x$ is a point where $u$ is infinitesimally generalized linear, we can also require that $|u(y_n)-u(x_n)| \leq c \epsilon \sd(x,x_n)$, where $c$ is a constant depending on $|\nabla u|(x)$. Hence,
        \[
        \frac{|f(i_g(y_n))-f(i_g(x_n))|}{\sd_g(i_g(y_n),i_g(x_n))} \geq \frac{|f(i_g(x))-f(i_g(x_n))|}{2c\epsilon \sd(x,x_n)} \geq \frac{|f(i_g(x))-f(i_g(x_n))|}{2c\epsilon \sd_g(i_g(x_n),i_g(x))}.
        \]
        Taking limits, we deuce $c \epsilon \ssf{L}(f) \geq \lip_g(f)(x)$. Since $\epsilon>0$ is arbitrary, we conclude.
    \end{proof}
\end{lemma}

The next theorem is the key step to prove Theorem \ref{T20}.

\begin{thm} \label{T15}
    Let $(\X,\sd,\m)$ be an $\RCD(K,N)$ space, let $\Omega \subset \X$ be open, and let $u \in \W^{1,1}(\Omega)$ be a solution of the minimal surface equation. Let $x \in \Omega$ be a point where $u$ is infinitesimally generalized linear and continuous, let $(\epsilon_n,\psi_n)$ be a blow-up of $\X$ at $x$, and let $(\epsilon_n,\psi_n,u^\infty)$ be the corresponding blow-up of $u$ given by Theorem \ref{T16}.
    Let $j:\bb{R}^k \to \Graph(u^{\infty}) \subset \bb{R}^k \times \bb{R}$ be the projection on the graph, i.e. $j(x):=(x,u^\infty(x))$.
    Then, there exists a sequence $\{\delta_n\}_{n \in \bb{N}}$ decreasing to zero such that the maps
    \[
    \psi'_n :=i_g \circ \psi_n \circ j^{-1}:j(\bar{B}^{\bb{R}^k}_1(0)) \to (\G(u) \cap \bar{B}^{\X}_{\epsilon_n}(x) \times \bb{R}, \epsilon_n^{-1}\sd_g)
    \]
    are $\delta_n$-$GH$ maps.
    \begin{proof}
    We first prove that for every $\delta>0$, if $n$ is sufficiently large, $\psi'_n (j(\bar{B}^{\bb{R}^k}_1(0)))$ is a $\delta$-net in $(\G(u) \cap \bar{B}^{\X}_{\epsilon_n}(x) \times \bb{R}, \epsilon_n^{-1}\sd_g)$.
    
    Assume by contradiction that this is not the case. Hence, there exist $\delta>0$ and a sequence of points $\bar{x}_n:=(x_n,t_n) \in \G(u) \cap \bar{B}^{\X}_{\epsilon_n}(x) \times \bb{R}$, such that $\epsilon_n^{-1}\sd_g(\bar{x}_n,\bar{y})>\delta$ for every $\bar{y}$ in the image of $i_g \circ \psi_n$.
    
    Consider the function $v:\Omega \to \bb{R}$ which coincides with $u$ outside of $\{x_n\}$ and such that $v(x_n)=t_n$. Note that the graph of this functions belongs to $\partial \Epi(u)$, so that we can apply Theorem \ref{T16}.
    Let $(\epsilon,\tilde{\psi}_n)$ be a blow-up of $\X$ at $x$ such that $\psi_n=\tilde{\psi}_n$ outside of a singleton $\{x'_n\}$ and such that $\tilde{\psi}_n(x'_n)=x_n$.

    By Theorem \ref{T16} (and Remark \ref{R2}), there exists a sequence $\delta_n \to 0$ such that, assuming for simplicity that $u(x)=0$, it holds
    \[
    |\epsilon_n^{-1} v(\tilde{\psi}_n(x'_n))-u^\infty(x'_n)| \leq \delta_n, \quad 
    |\epsilon_n^{-1} u(\psi_n(x'_n))-u^\infty(x'_n)| \leq \delta_n.
    \]
    Hence,
    \begin{equation} \label{E57}
    \epsilon_n^{-1}|t_n-u(\psi_n(x'_n))| \leq 2\delta_n.
    \end{equation}
    Since $\epsilon_n^{-1}\sd(x_n,\psi_n(x'_n)) \to 0$,
   it follows from \eqref{E57} that $\sd_g((x_n,t_n),i_g(\psi_n(x'_n)) \to 0$, a contradiction.
    \par 
    Now we need to show that $\psi'_n$ almost preserves distances. We pick two points $(x_1,u^{\infty}(x_1))$ and $(x_2,u^{\infty}(x_2))$ in the domain of $\psi'_n$ and we need to compute
    \[
    \epsilon_n^{-1}\sd_g (\psi'_n (x_1,u^{\infty} (x_1)),\psi'_n (x_2,u^{\infty}(x_2))).
    \]
    The previous expression is equal to
    \begin{align*}
    &\epsilon_n^{-1} \sd_\times ((\psi_n(x_1),u(\psi_n(x_1))),(\psi_n(x_2),u(\psi_n(x_2))))
    \\
    &=\sqrt{\epsilon_n^{-2} \sd(\psi_n(x_1),\psi_n(x_2))^2+\epsilon_n^{-2} |u(\psi_n(x_1))-u(\psi_n(x_2))|^2}.
    \end{align*}
    Using the inequality $|\sqrt{a^2+b^2}-\sqrt{c^2+d^2}| \leq |a-c|+|b-d|$, and denoting by $\sd_e$ Euclidean distances in any dimension, it holds
    \begin{align*}
    &|\epsilon_n^{-1}\sd_g (\psi'_n (x_1,u^{\infty} (x_1)),\psi'_n (x_2,u^{\infty}(x_2)))-\sd_e ((x_1,u^{\infty} (x_1)),(x_2,u^{\infty}(x_2)))|
    \\
    &\leq |\epsilon_n^{-1} \sd(\psi_n(x_1),\psi_n(x_2))-\sd_e(x_1,x_2)| +
    |\epsilon_n^{-1}u(\psi_n(x_1))-\epsilon_n^{-1}u(\psi_n(x_2))-u^{\infty}(x_1)+u^{\infty}(x_2)|.
    \end{align*}
    The last term goes uniformly to $0$ as $n$ increases, concluding the proof.
    \end{proof}
\end{thm}

Next we have an abstract lemma. Combining this with Theorem \ref{T15} and the previous results of this section, we then deduce Theorem \ref{T20}.

\begin{lemma} \label{LT16}
        Let $(\ssf{Y},\sd_y)$ be a metric space, let $y \in \ssf{Y}$, and let $f \in \Lip(\ssf{Y})$ be a function with
        \[
        \lip(f)(y)=\lim_{n \to + \infty}\frac{|f(y_n)-f(y)|}{\sd_y(y_n,y)}.
        \]
        Let $A \subset \bb{R}^k$ and $A_n \subset \ssf{Y}$ be sets such that $0 \in A $ and $y_n \in A_n$. Assume that $\epsilon_n,\delta_n$ are sequences decreasing to zero and that we have 
        $\delta_n$-$GH$ approximations $\psi_n:A \to (A_n,\epsilon_n^{-1} \sd_y)$ with $\psi_n(0)=y$ and $\psi_n(x_n)=y_n$ for some $x_n \in A$. If the sequence $x_n$ does not have $0$ as a limit point, and there exists a linear function $f^{\infty}: \bb{R}^k \to \bb{R}$ such that
        \[
        \|\epsilon_n^{-1}(f \circ \psi_n - f(y))-f^{\infty}\|_{\infty,A} \to 0, 
        \]
        then $\lip(f)(y)=\lip(f^{\infty})(0)$.
        \begin{proof}
            We first show that $\lip(f)(y) \leq \lip(f^{\infty})(0)$.
            Without loss of generality, we may suppose that $f(y)=0$, so that
            $
            |\epsilon_n^{-1}f(y_n)-f^{\infty}(x_n) | \to 0.
            $
            Note that $\epsilon_n^{-1} \sd(y,y_n)$ is bounded away from zero since $x_n$ has this property and $\psi_n$ are $\delta_n$-GH approximations. Hence, it holds
            \[
            \Big|\frac{f(y_n)}{\sd_y(y_n,y)}-\frac{f^{\infty}(x_n)}{\epsilon_n^{-1}\sd_y(y,y_n)} \Big| \to 0.
            \]
            Since $f^{\infty}$ is linear, to prove $\lip(f)(y) \leq \lip(f^{\infty})(0)$, it is sufficient to show that 
            \[
            \Big |\frac{f^{\infty}(x_n)}{\epsilon_n^{-1}\sd_y(y,y_n)} - \frac{f^{\infty}(x_n)}{|x_n|} \Big | \to 0.
            \]
            Since the sequence $x_n$ does not have $0$ as a limit point, the previous limit follows from the fact that the maps $\psi_n$ are $\delta_n$-GH approximations with $\psi_n(x_n)=y_n$.
            \par
            We now show that $\lip(f^{\infty})(0) \leq \lip(f)(y)$. Since $f^{\infty}$ is linear, there exists $x_{\infty} \in \bb{R}^k$ such that $|x_{\infty}|=1$ and $\lip(f^{\infty})(0)=f^{\infty}(x_{\infty})$.
            By hypothesis, assuming as before that $f(y)=0$, it holds
            $
            |{\epsilon_n}^{-1}f(\psi_n(x_{\infty}))-f^{\infty}(x_{\infty}) | \to 0,
            $
            and $\sd_y(\psi_n(x_{\infty}),y) \leq \epsilon_n(1+\delta_n)$. Therefore
            \[
            \lip(f^{\infty})(0)=f^{\infty}(x_{\infty})=
            \lim_{n \to + \infty} {\epsilon_n}^{-1}|f(\psi_n(x_{\infty}))|
            \]
            \[
            \leq
            \limsup_{n \to + \infty} \frac{(1+\delta_n)|f(\psi_n(x_{\infty}))|}{\sd_y(\psi_n(x_{\infty}),y)} \leq \lip(f)(y),
            \]
            concluding the proof.
        \end{proof}
\end{lemma}

The next result is the key theorem of the section, which then implies Theorems \ref{CT5} and \ref{CT6} from the Introduction.

\begin{thm} \label{T20}
    Let $(\X,\sd,\m)$ be an $\RCD(K,N)$ space, let $\Omega \subset \X$ be an open set, let $u \in \W^{1,1}(\Omega)$ be a solution of the minimal surface equation, and let $\theta \in \Lip(\G(u))$. Then, $\theta \circ i_g \in \W^{1,1}(\Omega)$ and, for $\m$-almost every $x \in \Omega$, it holds
    \[
    \lip_g(\theta)^2 \circ i_g=|\nabla \theta \circ i_g|^2-\frac{(\nabla \theta \circ i_g \cdot \nabla u)^2}{W_u^2}.
    \]
    \begin{proof}
        By Lemma \ref{L15}, it holds $\theta \circ i_g \in \W^{1,1}(\Omega)$.
        By Theorem \ref{T14}, for every $\epsilon>0$, there exists $f \in \Lip(\Omega)$, such that $\m\{f \neq \theta \circ i_g\} \leq \epsilon$. 
        Let $x \in \Omega$ satisfy the following.
        \begin{enumerate}
            \item $\{f = \theta \circ i_g\}$ has density $1$ in $x$. 
            \item \label{item1} $f$ is infinitesimally generalized linear at $x$.
            \item $u$ satisfies the property of Theorem \ref{T16} in $x$ (in particular, it is is infinitesimally generalized linear at $x$ as well).
            \item $\lip(f)(x)=|\nabla f|(x)=|\nabla (\theta \circ i_g)|(x)$.
            \item $\nabla \theta \circ i_g \cdot \nabla u=\nabla f \cdot \nabla u$.
        \end{enumerate}
        Consider the function $f \circ i_g^{-1}:\Graph(u) \to \bb{R}$. It is easy to check that $f \circ i_g^{-1}$ is Lipschitz, so that we consider its Lipschitz extension $f_g:\G(u) \to \bb{R}$.
        By Lemma \ref{L16}, it holds 
        \begin{equation} \label{E55}
        \lip_g(f_g)(i_g(x))=\lip_g(\theta)(i_g(x)).
        \end{equation}
        By Lemma \ref{L17}, there exist points $x_n \to x$ such that
        \[
        \lip_g(f_g)(i_g(x))=\lim_{n \to + \infty} \frac{|f_g(i_g(x_n))-f_g(i_g(x))|}{\sd_g(i_g(x_n),i_g(x))}.
        \]
        Let $\epsilon_n:=\sd(x,x_n)$ and, modulo passing to a subsequence, let $(\epsilon_n,\psi_n)$ be a blow-up of $\X$ in $x$ such that $x_n \in \ssf{Im}(\psi_n)$.
        Let $f^\infty$ and $u^\infty$ be the blow-ups of $f$ and $u$ in $x$ along (a subsequence of) $(\epsilon_n,\psi_n)$. Let $j:\bb{R}^k \to \Graph(u^{\infty}) \subset \bb{R}^k \times \bb{R}$ be the projection on the graph of $u^\infty$, i.e. $j(x):=(x,u^\infty(x))$. 
        By Theorem \ref{T15}, the maps
        \[
    \psi'_n :=i_g \circ \psi_n \circ j^{-1}:j(\bar{B}^{\bb{R}^k}_1(0)) \to (\G(u) \cap \bar{B}^{\X}_{\epsilon_n}(x) \times \bb{R}, \epsilon_n^{-1}\sd_g)
    \]
    are $\delta_n$-Gromov Hausdorff approximations for some sequence $\delta_n \to 0$. 
    Moreover, by definition of $f^\infty$, it holds
    \[
    \|\epsilon_n^{-1}( f_g \circ i_g \circ \psi_n \circ j^{-1}-f_g(i_g(x)))-f ^{\infty} \circ j^{-1} \| _{\infty,j(\bar{B}^{\bb{R}^k}_1(0))}
    \]
    \[
    = 
    \|\epsilon_n^{-1}( f \circ \psi_n  -f(x))-f ^{\infty}  \| _{\infty,\bar{B}^{\bb{R}^k}_1(0)}
    \to 0.
    \]
    In addition, since $u^\infty$ is linear, by our choice of $\epsilon_n$, the set $j(\psi_n^{-1}(x_n)) \subset \Graph(u^\infty)$ cannot have $0$ as a limit point.
    Hence, using Lemma \ref{LT16}, it holds
    \[
    \lip_g(f_g)^2(i_g(x))
    =\lip(f ^{\infty} \circ j^{-1})^2.
    \]
    By a computation in $\bb{R}^k$, it follows
    \[
    \lip_g(f_g)^2(i_g(x))
    =\lip(f ^{\infty} \circ j^{-1})^2
    ={\lip}(f^{\infty} )^2-\frac{(\lip(f ^{\infty}) \cdot \lip(u^{\infty}))^2}{\sqrt{1+\lip(u^{\infty})^2}}.
    \]
    Since $x$ is a point where $f$ and $u$ are infinitesimally generalized linear, by Corollary \ref{CL19}, it holds $\nabla f^\infty \cdot \nabla u^\infty=\nabla f \cdot \nabla u=\nabla \theta \circ i_g \cdot \nabla u$.
    Combining with \eqref{E55}, it holds
    \[
    \lip_g(\theta)^2 (i_g(x))=
    \Big(|\nabla \theta \circ i_g |^2-\frac{(\nabla \theta \circ i_g \cdot \nabla u)^2}{W_u^2}\Big)(x),
    \]
    concluding the proof. 
    \end{proof}
\end{thm}

\begin{corollary} \label{CP7}
     Let $(\X,\sd,\m)$ be an $\RCD(K,N)$ space, let $\Omega \subset \X$ be open, and let $u \in \W^{1,1}_{loc}(\Omega)$ be a solution of the minimal surface equation. Let $\theta_1, \theta_2 \in \Lip(\G(u))$.
    For $\m$-a.e. $x \in \Omega$, it holds
    \[
    \Big( \lip_g(\theta_1 ) \cdot \lip_g(\theta_2 )
    \Big) \circ i_g
    = \nabla (\theta_1 \circ i_g) \cdot \nabla (\theta_2 \circ i_g)-\frac{1}{W_u^2}(\nabla (\theta_1 \circ i_g) \cdot \nabla u)(\nabla (\theta_2 \circ i_g)\cdot \nabla u).
    \]
    \begin{proof}
        It follows immediately from the previous theorem.
    \end{proof}
\end{corollary}

We now recall the definition of a product of gradients in a generic metric space. This will allow to consider products of gradients on a minimal graph.

\begin{definition}
    Let $(\ssf{Y},\sd)$ be a metric space and $f,g:\ssf{Y} \to \bb{R}$ be functions with finite local Lipschitz constant in $y \in \ssf{Y}$. We define
    \[
    \lip(f) \cdot \lip(g):=\frac{1}{4}(\lip(f+g)^2-\lip(f-g)^2).
    \]
\end{definition}

Observe that the object defined in the previous definition may fail to be a quadratic form on generic metric spaces. This additional regularity, in our setting, is a consequence of Corollary \ref{CP7}.

\begin{customthm}{5}
    Let $(\X,\sd,\m)$ be an $\RCD(K,N)$ space, let $\Omega \subset \X$ be open, and let $u \in \W^{1,1}_{loc}(\Omega)$ be a solution of the minimal surface equation. The function $\cdot:\Lip(\G(u)) \times \Lip(\G(u)) \to \sL^\infty(\G(u))$ given by 
    \[
    (\phi_1,\phi_2) \mapsto \lip_g(\phi_1) \cdot \lip_g(\phi_2)
    \]
    is symmetric, bilinear, it satisfies the chain rule and the Leibniz rule in both entries and $\lip_g (\phi_1) \cdot \lip_g(\phi_2) \leq \lip_g(\phi_1) \lip_g(\phi_2)$.
    \begin{proof}
        The fact that $\lip_g(\phi_1) \cdot \lip_g(\phi_2) \leq \lip_g(\phi_1) \lip_g(\phi_2)$ follows from the fact that $\lip_g(\phi_1+\phi_2) \leq 
        \lip_g(\phi_1)+\lip_g(\phi_2)$ and $\lip_g(\phi_1) \leq 
        \lip_g(\phi_1-\phi_2)+\lip_g(\phi_2)$. All the other properties follow from the representation given by Corollary \ref{CP7}.
    \end{proof}
\end{customthm}

Let $(\X,\sd,\m)$ be an $\RCD(K,N)$ space, let $\Omega \subset \X$ be open, and let $u \in \W^{1,1}_{loc}(\Omega)$ be a solution of the minimal surface equation. We set $u_g:\G(u) \to \bb{R}$ to be the height function 
\begin{equation} \label{E72}
    u_g(x,t):=t,
\end{equation}
    so that $u_g \circ i_g=u$. 



\begin{customthm}{6}
Let $(\X,\sd,\m)$ be an $\RCD(K,N)$ space, let $\Omega \subset \X$ be open, and let $u \in \W^{1,1}_{loc}(\Omega)$ be a solution of the minimal surface equation.
Let $\phi \in \Lip_c(\G(u))$, then ${\lip}_g(\phi ) \cdot {\lip}_g(u_g) \in \sL^1(\G(u))$, and
\[
\int_{\G(u)} \lip_g(\phi) \cdot \lip_g(u_g) \, d\m_g=0.
\]
\begin{proof}
    Since $\phi \in \Lip_c(\G(u))$ and $u_g \in \Lip(\G(u))$, then ${\lip}_g(\phi) \cdot {\lip}_g(u_g) \in \sL^1(\G(u))$ trivially.
    Since $\m_g(B \times \bb{R})=\int_B W_u \, d\m$ for every Borel set $B \subset \Omega$, it holds
    \[
    \int_{\G(u)} \lip_g(\phi) \cdot \lip_g(u_g) \, d\m_g
    =
    \int_{\Omega} \Big( \lip_g(\phi) \cdot \lip_g(u_g) \Big) \circ i_g \, W_u d \m.
    \]
    By Corollary \ref{CP7}, it follows
    \begin{equation} \label{E62}
    \int_{\G(u)} \lip_g(\phi) \cdot \lip_g(u_g) \, d\m_g
    =
    \int_{\Omega} \frac{\nabla (\phi \circ i_g) \cdot \nabla u}{W_u} d \m.
    \end{equation}
    The function $\phi \circ i_g$ has compact support in $\Omega$ and, by Lemma \ref{L15}, it holds $\phi \circ i_g \in \W^{1,1}(\Omega)$. Since $u$ solves the minimal surface equation in $\Omega$, the right hand side of \eqref{E62} vanishes, concluding the proof.
\end{proof}
\end{customthm}

Using the chain rule, the Leibniz rule, and applying the previous theorem, we also obtain the following corollary. This is the version of 'integration by parts' that we will use repeatedly in the next section.
\begin{corollary} \label{CCC10} 
Let $(\X,\sd,\m)$ be an $\RCD(K,N)$ space, let $\Omega \subset \X$ be open, and let $u \in \W^{1,1}_{loc}(\Omega)$ be a solution of the minimal surface equation.
    Let $\phi \in \Lip_c(\G(u))$, let $A \subset \bb{R}$ be an open set containing the image of $u_g$, and let $h=g \circ u_g$ with $g: A \to \bb{R}$ smooth. Then
    \[
    \int_{\G(u)} \lip_g(\phi) \cdot \lip_g(h) \, d\m_g =- \int_{\G(u)} \lip_g(u_g )^2 g''(u_g) \phi\, d\m_g.
    \]
    \begin{proof}
        By the chain rule and the Leibniz rule, it holds
        \begin{align*}
            & \int_{\G(u)} \lip_g(\phi) \cdot \lip_g(h) \, d\m_g =
            \int_{\G(u)} g'(u_g)\lip_g(\phi) \cdot \lip_g(u_g) \, d\m_g  \\
            &
            =\int_{\G(u)} \lip_g(g'(u_g)\phi) \cdot \lip_g(u_g) \, d\m_g -
            \int_{\G(u)} \lip_g(u_g )^2 g''(u_g) \phi\, d\m_g.
        \end{align*}
        The first summand in the previous expression vanishes by Theorem \ref{CT6}, concluding the proof.
    \end{proof}
\end{corollary}

\section{Bernstein Property} \label{S3}
The first subsection contains the proof of Theorem \ref{CT1} (in a stronger version, see Theorem \ref{T13}), while the second one contains the proofs of Theorems \ref{CT2} and \ref{CT3} (as anticipated we prove the stronger version of Theorem \ref{CT3} given by Corollary \ref{C12}).

\subsection{Harnack inequality for the height function on a minimal graph}

Let $(\X,\sd,\m)$ be an $\RCD(0,N)$ space, let $B_R(p) \subset \X$ be a ball, and let $u \in \W^{1,1}_{loc}(B_R(p))$ be a solution of the minimal surface equation.
The goal is to mimic the strategy used in \cite{Ding} to prove the Harnack inequality for $u_g$ on $\G(u)$ i.e. Theorem \ref{T10}. 
The challenge in adapting the aforementioned strategy is to prove that $u_g$ is harmonic on $\G(u)$ and this was done in Section \ref{S2}. 
Besides harmonicity on the graph, the other key steps are the validity of Poincaré and Sobolev inequalities on $\G(u)$ i.e. Theorems \ref{T8} and \ref{T9}. \par 
These theorems can be obtained in our setting with the same ideas used in \cite{Ding}. For this reason we only give a detailed proof of Theorem \ref{T8}, so that one sees what changes need to be made from the corresponding proof in \cite{Ding}; the same exact changes allow then to obtain also the Sobolev inequality from the analogous result in the aforementioned work.
Once the harmonicity of $u$ on $\G(u)$ and the Poincaré and Sobolev inequalities are established, the Harnack inequality follows formally repeating the same argument of \cite{Ding} and for this reason the proof of
Theorem \ref{T10} is only sketched. 

\begin{definition} \label{D6}
Let $(\X,\sd,\m)$ be an $\RCD(K,N)$ space, let $\Omega \subset \X$ be open, and let $u \in \W^{1,1}_{loc}(\Omega)$ be a solution of the minimal surface equation.
Let $ (x,t), (y,s) \in \Omega \times \bb{R}$ and define $\rho_{(x,t)}(y,s):=\sd(x,y)+|t-s|$. 
Given $\bar{x} \in \Omega \times \bb{R}$, we set 
\[
D_{\bar{x},r}:=\{\bar{y} \in \Omega \times \bb{R}: \rho_{\bar{x}}(\bar{y}) < r \}
\]
and 
\[
A_{\bar{x},r}:=\{x \in \Omega: \{x\} \times \bb{R} \cap \G(u) \cap D_{\bar{x},r} \neq \emptyset\}.
\]
For every $t \in \bb{R}$, we also set 
    \[
    E_{\bar{x},t}:=\{\bar{y} \in \G(u):u_g(\bar{y})>t+u_g(\bar{x})\}
    \]
    and
    \[
    E'_{\bar{x},t}:=\{\bar{y} \in \G(u):u_g(\bar{y}) < t+u_g(\bar{x})\}.
    \]
\end{definition}

\begin{lemma} \label{L10}
Let $(\X,\sd,\m)$ be an $\RCD(0,N)$ space, let $B_R(p) \subset \X$ be a ball, and let $u \in \W^{1,1}_{loc}(B_R(p))$ be a solution of the minimal surface equation. There exists a constant $C(N)>0$ such that, for every $(x,t_x)=:\bar{x} \in \G(u)$, and
$r,s \in (0,R-\sd(x,p))$ with $r>s$, it holds
    \[
    \m_g(D_{\bar{x},s}) \geq  C\m_g(D_{\bar{x},r})\frac{s^N}{r^N}
    \]
    and
    \begin{equation} \label{E63}
    \frac{1}{C}\m(B_s(x)) \leq \m_g(D_{\bar{x},s})  \leq C\m(B_s(x)).
    \end{equation}
    \begin{proof}
    We denote by $c_1,c_2 \cdots$ positive constants depending on $N$.
        From Theorem \ref{T7}, the fact that $\m_\times$ is doubling and the fact that the ball $B_s(\bar{x}) \subset \X \times \bb{R}$ contains the rectangle $B_{s/2}(x) \times (t_x-s/2,t_x+s/2)$, it holds
        \begin{align*}
        \m_g(D_{\bar{x},s}) \geq \m_g(B_{s/2}(\bar{x}))
        \geq c_1 \frac{2\m_\times(B_{s/2}(\bar{x}))}{s}
        \geq c_2 \frac{\m_\times(B_{s}(\bar{x}))}{s}
        \\
         \geq c_2 \frac{\m(B_{s/2}(x))s/2}{s}
        \geq c_3 \m(B_s(x)).
        \end{align*}
        By a similar argument, one obtains
        \[
        \m(B_s(x))
        \geq c_4\m(B_{2s}(x))
        \geq
        c_5 \frac{\m_\times(B_{s}(\bar{x}))}{s}
        \geq
        c_6 \m_g(D_{\bar{x},s}),
        \]
        proving \eqref{E63}.
        Since $\m$ is doubling, by a standard argument (see \cite[Equation $(2.1)$]{Ding}),
        it holds 
        \[
        c_3 \m(B_s(x)) \geq 
        c_7\m(B_r(x))\frac{s^N}{r^N}.
        \]
        Combining with \eqref{E63}, the statement follows.
    \end{proof}
\end{lemma}

\begin{lemma} \label{L14}
    Let $(\X,\sd,\m)$ be an $\RCD(0,N)$ space, let $B_R(p) \subset \X$ be a ball, and let $u \in \W^{1,1}_{loc}(B_R(p))$ be a solution of the minimal surface equation. There exists a constant $C(N)>0$ such that, for every $(x,t_x)=:\bar{x} \in \G(u)$ and every $r>0$ with $B_{3r}(x) \subset B_R(p)$, for $\lambda^1$-a.e. $ t \in (-r,r)$ it holds
    \begin{equation} \label{E34}
             \m(\{u>u_g(\bar{x})+t\} \cap B_{r-t}(x)) \geq C \, \m_g(E_{\bar{x},t} \cap D_{\bar{x},r}) .
        \end{equation}
    Moreover, if 
    \begin{equation} \label{E35}
    \m_g(E_{\bar{x},t} \cap D_{\bar{x},r}) \leq \frac{1}{2} \m_g( D_{\bar{x},r}),
    \end{equation}
    then
    \begin{equation} \label{E36}
            \m(\{u < u_g(\bar{x})+t\} \cap B_{r+|t|}(x)) \geq C \, \m(B_r(x)).
        \end{equation}
    Analogous inequalities hold replacing $E_{\bar{x},t}$ with $E_{\bar{x},t}'$ and reversing the inequality signs appearing in the left hand sides of inequalities \eqref{E34} and \eqref{E36}.
    \begin{proof}
        We first prove \eqref{E36}. Modulo requiring that $t$ is outside of a negligible set, we can assume that $P(\Epi(u),\X \times \{t+t_x\}=0)$.
        We denote $\bar{x}_t:=(x,t_x+t)$ and we consider the compact set $V \subset \X \times \bb{R}$ given by the closure of
        \[
        \Epi(u) \cap D_{\bar{x}_t,r+|t|} \cap \X \times (-\infty,t_x+t)
        \]
        and we define the competitor $C:=\Epi'(u) \cup V$. 
        
\begin{figure}[h]
\includegraphics[width=7cm, height=5cm]{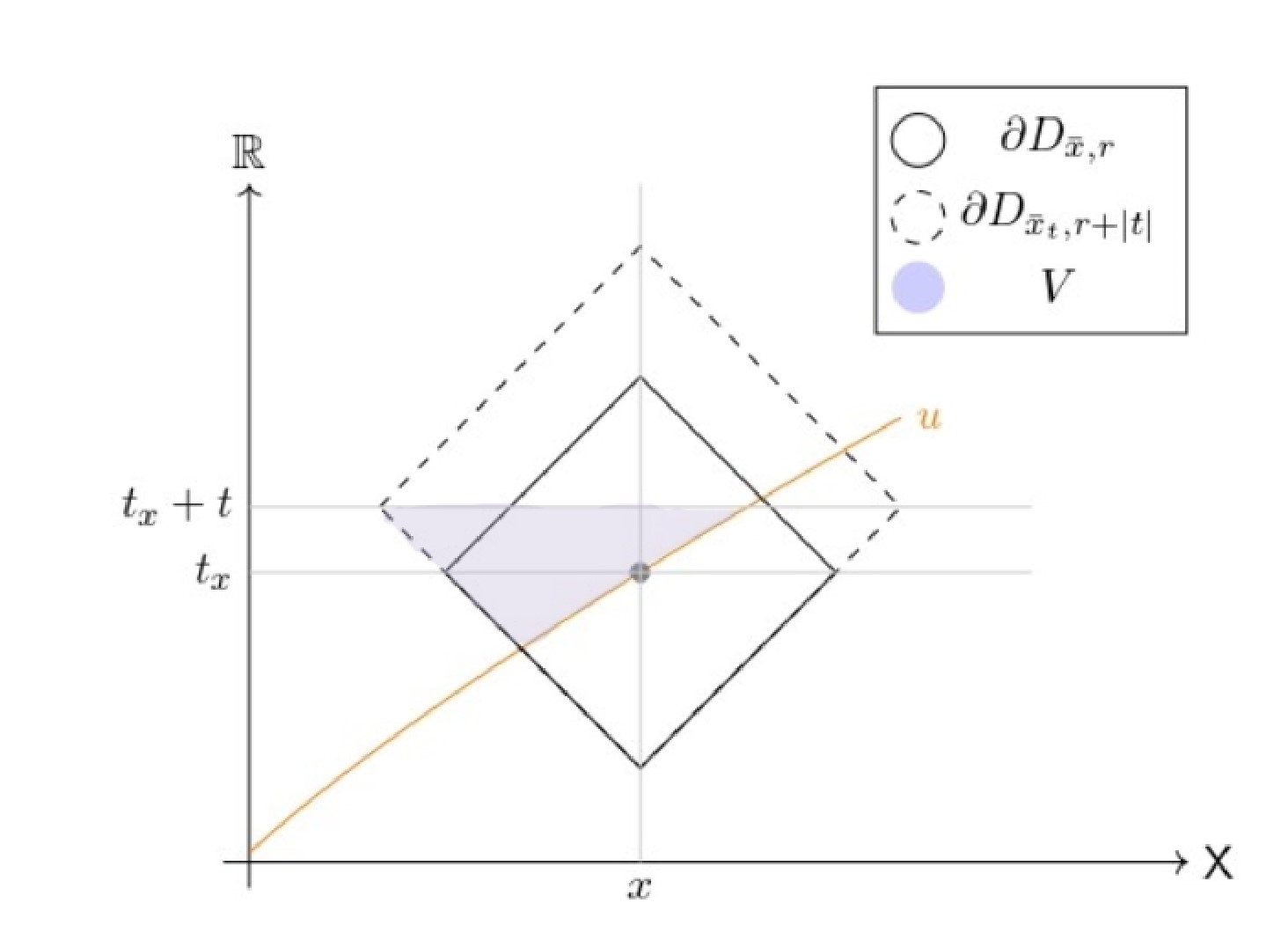}
\centering
\end{figure}

Since $\Epi'(u)$ minimizes the perimeter, it holds
        \begin{equation} \label{E37}
        P(\Epi'(u), D_{\bar{x},3r}) \leq P(C,D_{\bar{x},3r}).
        \end{equation}
        We now claim that 
        \begin{equation} \label{E38}
        P(\Epi'(u), (\G(u) \setminus \partial V) \cap D_{\bar{x},3r})
        =
        P(C, (\partial C \setminus \partial V) \cap D_{\bar{x},3r}).
        \end{equation}
        To prove the claim observe that, by the definition of $V$, we have $\partial C \setminus \partial V=\G(u) \setminus \partial V$ and, if $q \in \G(u) \setminus \partial V$, then in a small ball centered at $q$ we have that $C=\Epi'(u)$, which implies the claim. \par 
        Subtracting \eqref{E38} from \eqref{E37}, we deduce that
        \begin{equation} \label{E64}
        P(\Epi'(u),\partial V \cap \G(u) \cap D_{\bar{x},3r}) \leq P(C, \partial V \cap \partial C \cap D_{\bar{x},3r}).
        \end{equation}
        We denote by $c_1.c_2, \cdots $ constants depending only on $N$. Using first the definition of $V$, then \eqref{E35}, and finally Lemma \ref{L10}, we obtain
        \begin{align*}
        P(\Epi'(u),\partial V \cap \G(u) \cap D_{\bar{x},3r}) \geq \m_g(E'_{\bar{x},t} \cap D_{\bar{x},r})
        \\
        \geq c_1 \m_g( D_{\bar{x},r})
        \geq c_2 \m(B_r(x)),
        \end{align*}
        so that, combining with \eqref{E64}, to prove \eqref{E36} it is enough to show that
        \[
        P(C, \partial V \cap \partial C \cap D_{\bar{x},3r}) \leq c_3 \m(\{u < u_g(\bar{x})+t\} \cap B_{r+|t|}(x)).
        \]
        To this aim observe that $\partial V \cap \partial C$ is the disjoint union of
        \[
        A:=\X \times \{t_x+t\} \cap \Epi(u) \cap D_{\bar{x}_t,r+|t|},
        \]
        \[
        B:=\X \times (-\infty,t_x+t) \cap \Epi(u) \cap \partial D_{\bar{x}_t,r+|t|},
        \]
        and
        \[
        D:=(\bar{A}\setminus A) \cup (\bar{B}\setminus B).
        \]
        If $q \in A$, then in a sufficiently small ball centered at $q$ the set $C$ coincides with the subgraph of the constant function $t+t_x$. This, together with the area formula given by Proposition \ref{P1}, implies that 
        \[
        P(C,A) \leq \m(\{u < t_x+t\} \cap B_{r+|t|}(x)).
        \]
        Similarly, if $q \in B$, then in a sufficiently small ball centered at $q$ the set $C$ coincides with the epigraph of a function with Lipschitz constant $1$. Moreover if $q \in B$, then its projection on $\X \times \{t+t_x\}$ belongs to $A$ (this is clear from the picture, and depends on the fact that we defined $V$ using $D_{\bar{x}_t,r+|t|}$ and not $D_{\bar{x},r+|t|}$). This implies that $P(C,B) \leq c_4 P(C,A)$. 
        \par 
        Finally for $\lambda^1$-a.e. $t \in (-r,r)$ we have that $\m_g(\X \times \{t_x+t\})=0$ and
        $\m_g(\partial D_{\bar{x},r+|t|})=0$. Using the area formula given by Proposition \ref{P1} we get that for every such $t$, it holds $\m(\pi_\X(\bar{A} \setminus A))=\m(\pi_\X(\bar{B} \setminus B))=0$, which in turn implies that $P(C,D)=0$. \par
        In particular,
        \[
        P(C,\partial V \cap \partial C)=P(C,A \cup B) \leq c_5 \m(\{u < t_x+t\} \cap B_{r+|t|}(x)),
        \]
        as desired. \par 
        The proof of the other identity follows an identical argument replacing $V$ with the closure of
        \[
        \Epi'(u) \cap D_{\bar{x}_t,r-t} \cap \X \times (t_x+t,\infty).
        \]
        In this case, since since the r.h.s. is $\m_g(E_{\bar{x},t})$, we do not need to use the condition \eqref{E35}, which was needed to say that $\m_g(E'_{\bar{x},t}) \geq c_1\m(B_r(x))$.
    \end{proof}
    \end{lemma}

The next lemma corresponds to Lemma $3.3$ in \cite{Ding}. Our formulation is slightly different and, later on, it will allow us to avoid the use of the Coarea formula on $\G(u)$, since this tool is a priori not available.
\begin{lemma} \label{L9}
Let $(\X,\sd,\m)$ be an $\RCD(0,N)$ space, let $B_R(p) \subset \X$ be a ball, and let $u \in \W^{1,1}_{loc}(B_R(p))$ be a solution of the minimal surface equation.
    There exists a constant $C(N)>0$ such that for every $(x,t_x)=:\bar{x} \in \G(u)$ and $r>0$ with $B_{3r}(x) \subset B_R(p)$ the following happens. For $\lambda^1$-a.e. $ t \in (-r,r)$,
    if 
    \[
    \m_g(E_{\bar{x},t} \cap D_{\bar{x},r}) \leq \frac{1}{2} \m_g( D_{\bar{x},r}),
    \]
    then
    \[
    \m_g(E_{\bar{x},t} \cap D_{\bar{x},r}) \leq Cr P(\{u>u_g(\bar{x})+t\} , A_{\bar{x},3r}).
    \]
    The same statement holds replacing $E_{\bar{x},t}$ with $E_{\bar{x},t}'$ and reversing the inequality sign appearing in the right hand side of the last inequality.
    \begin{proof}
        We denote by $c_1,c_2 \cdots $ quantities depending only on $N$. By the choice of $t$ and $r$, it holds
        \[
        3rP(\{ u>t_x+t\} , B_{r+|t|}(x)) \geq
        (r+|t|)P(\{u>t_x+t\} , B_{r+|t|}(x)).
        \]
        Thanks to Proposition \ref{P9}, this last quantity is greater than or equal to
        \begin{equation} \label{E51}
         c_1\min{\{\m(\{u>t_x+t\} \cap B_{r+|t|}(x)),\m(\{u \leq t_x+t\} \cap B_{r+|t|}(x) )\}}.
        \end{equation}
        Because of \eqref{E36} in Lemma \ref{L14}, for $\lambda^1$-a.e. $t \in (-r,r)$ the quantity in \eqref{E51} is greater than
        $
        c_2 \m(\{u > t+t_x \} \cap B_{r+|t|}(x))
        $, which is then greater than
        $
        c_3\m_g(E_{\bar{x},t} \cap D_{\bar{x},r})
        $ by equation \eqref{E34} in the same lemma. 
        Hence, we proved that
        \[
        \m_g(E_{\bar{x},t} \cap D_{\bar{x},r}) \leq c_4r P(\{u>t_x+t\} , B_{r+|t|}(x)),
        \]
        so that to conclude we only need to show that
        \[
        P(\{u>t_x+t\} , B_{r+|t|}(x)) \leq  P(\{u>t_x+t\} , A_{\bar{x},3r}).
        \]
        To this aim we will show that $\partial \{u>t_x+t\} \cap B_{r+|t|}(x) \subset A_{\bar{x},3r}$, from which the desired inequality will follow immediately. Let $y \in \partial \{u>t_x+t\} \cap B_{r+|t|}(x)$ and assume by contradiction that $(y,t_x+t) \notin \G(u)$. As $\G(u)$ is closed, there exists a small ball centered in $(y,t_x+t)$ which is either fully contained in the interior of $\Epi(u)$ or in the interior of its complement. We only consider the first case, the other one being analogous. Call $B$ the aforementioned ball, and let $\epsilon >0$ be small enough such that $B_\epsilon (y) \times [t_x+t-\epsilon,t_x+t+\epsilon] \subset B$, so that we deduce that $u<t_x+t-\epsilon$ in $B_\epsilon (y)$, contradicting the fact that $y \in \partial \{u>t_x+t\}$. \par 
        We deduced that $(y,t_x+t) \in \G(u)$, so that to conclude we only need to note that $(y,t_x+t) \in D_{\bar{x},3r}$ since $t \in (-r,r)$.
    \end{proof}
\end{lemma}

\begin{lemma} \label{L12}
    Let $(\X,\sd,\m)$ be an $\RCD(0,N)$ space, let $B_R(p) \subset \X$ be a ball, and let $u \in \W^{1,1}_{loc}(B_R(p))$ be a solution of the minimal surface equation. Let $f:\bb{R} \to \bb{R}$ be a smooth function and define $\phi:\G(u) \to \bb{R}$ as $\phi:=f \circ u_g$. Let $(x,t_x)=:\bar{x} \in \G(u)$ and $r>0$ be such that $B_r(x) \subset B_R(p)$. Then
        \[
        \int_{ A_{\bar{x},r}} |\nabla (\phi \circ i_g)| \, d\m
        = \int_{ D_{\bar{x},r}} \lip_g( \phi) \, d\m_g.
        \]
        \begin{proof}
        By Theorem \ref{T20}, it holds $\lip_g(\phi) \circ i_g=|\nabla \phi \circ i_g|/W_u$ for $\m$-a.e. $x \in B_R(p)$. Since $\m_g(B \times \bb{R})=\int_B W_u \, d\m$ for every Borel set $B \subset B_R(p)$, and $D_{\bar{x},r} \cap \G(u)=A_{\bar{x},r} \times \bb{R} \cap \G(u)$, we obtain
        \[
         \int_{ D_{\bar{x},r}} \lip_g( \phi) \, d\m_g=
         \int_{ A_{\bar{x},r} \times \bb{R}} \lip_g( \phi) \, d\m_g
         =
         \int_{ A_{\bar{x},r}} |\nabla (\phi \circ i_g)| \, d\m,
        \]
        concluding the proof.
        \end{proof}
\end{lemma}

Given a function $f \in \sL^1(\G(u))$, we use the notation
\[
f_{\bar{x},r}:=\fint_{D_{\bar{x},r}}f \, d\m_g.
\]
\begin{thm} \label{T8}
Let $(\X,\sd,\m)$ be an $\RCD(0,N)$ space, let $B_R(p) \subset \X$ be a ball, and let $u \in \W^{1,1}_{loc}(B_R(p))$ be a solution of the minimal surface equation. There exists a constant $C(N)>0$ such that, for every $(x,t_x)=:\bar{x} \in \G(u)$, $r>0$ with $B_{3r}(x) \subset B_R(p)$, and every smooth monotone function $f:\bb{R} \to \bb{R}$, defining $\phi:\G(u) \to \bb{R}$ as $\phi:=f \circ u_g$, it holds
    \[
    \int_{D_{\bar{x},r}}|\phi -\phi_{\bar{x},r}| \, d\m_g \leq Cr \int_{ D_{\bar{x},3r}} \lip_g (\phi) \, d\m_g.
    \]
    \begin{proof}
    Observe that since $f$ is smooth and monotone, for every $s \in \bb{R}$, it holds $\{\phi >s\}=\{u_g > f^{-1}(s)\}$ (or $\{u_g < f^{-1}(s)\}$ if $f$ is decreasing) and if $A \subset \bb{R}$ has full $\lambda^1$ measure then also $f^{-1}(A)$ has this property. This will allow us to use Lemma \ref{L9} in what follows. We assume for simplicity that $f$ is increasing, as the other case is analogous. \par
    Suppose again for simplicity (the other case being again analogous) that 
    \[
    \m_g(\{\phi > \phi_{\bar{x},r}\} \cap D_{\bar{x},r}) \leq \m_g(\{\phi < \phi_{\bar{x},r}\} \cap D_{\bar{x},r}).
    \]
    In this case, for every $t \geq 0$, it holds
    \begin{align*}
    & \m_g(\{\phi > \phi_{\bar{x},r}+t\} \cap D_{\bar{x},r}) \leq 
    \m_g(\{\phi > \phi_{\bar{x},r}\} \cap D_{\bar{x},r})
    \\
    & \leq 
    \m_g(\{\phi < \phi_{\bar{x},r}\} \cap D_{\bar{x},r})
    \leq
    \m_g(\{\phi < \phi_{\bar{x},r}+t\} \cap D_{\bar{x},r}),
    \end{align*}
    so that in particular
    \[
    \m_g(\{\phi > \phi_{\bar{x},r}+t\} \cap D_{\bar{x},r}) \leq \frac{1}{2}\m_g(D_{\bar{x},r}).
    \]
    Hence, using first a variation of Fubini's Theorem (see \cite[Proposition $1.78$]{AmbFuscPall}) and then Lemma \ref{L9}, we obtain
        \begin{align*}
        \int_{\{\phi > \phi_{\bar{x},r}\} \cap  D_{\bar{x},r}}(\phi -\phi_{\bar{x},r}) \, d\m_g=\int_0^{+\infty} \m_g(\{\phi > t+\phi_{\bar{x},r}\} \cap D_{\bar{x},r}) \, dt,
        \\
        \leq Cr \int_0^{+\infty} P(\{\phi \circ i_g >t+ \phi_{\bar{x},r}\} , A_{\bar{x},3r})\, dt.
        \end{align*}
        This last expression, by the Coarea formula and Lemma \ref{L12}, is less than 
        \[
        Cr \int_{ A_{\bar{x},3r}} |\nabla (\phi \circ i_g)| \, d\m
        = Cr\int_{ D_{\bar{x},3r}} \lip_g( \phi) \, d\m_g.
        \]
        Summing up, we proved that
         \[
        \int_{\{\phi > \phi_{\bar{x},r}\} \cap  D_{\bar{x},r}}(\phi -\phi_{\bar{x},r}) \, d\m_g
        \leq Cr \int_{ D_{\bar{x},3r}} \lip_g( \phi) \, d\m_g.
        \]
        To conclude, we observe that
        \begin{align*}
        \int_{  D_{\bar{x},r}} &|\phi -\phi_{\bar{x},r}| \, d\m_g
        \\
        &=
        \int_{\{\phi > \phi_{\bar{x},r}\} \cap  D_{\bar{x},r}}(\phi -\phi_{\bar{x},r}) \, d\m_g
        +
        \int_{\{\phi < \phi_{\bar{x},r}\} \cap  D_{\bar{x},r}}(\phi_{\bar{x},r}-\phi) \, d\m_g
        \\
        &=2\int_{\{\phi > \phi_{\bar{x},r}\} \cap  D_{\bar{x},r}}(\phi -\phi_{\bar{x},r}) \, d\m_g
        \leq 2Cr \int_{ D_{\bar{x},3r}} \lip_g( \phi) \, d\m_g.
        \end{align*}
    \end{proof}
\end{thm}

We now state the Sobolev isoperimetric inequality.

\begin{thm} \label{T9}
Let $(\X,\sd,\m)$ be an $\RCD(0,N)$ space, let $B_R(p) \subset \X$ be a ball, and let $u \in \W^{1,1}_{loc}(B_R(p))$ be a strictly positive solution of the minimal surface equation. There exists a constant $C(N)>0$ satisfying the following. Let $(x,t_x)=:\bar{x} \in \G(u)$, $r \geq \tau >0$ with $B_{2r}(x) \subset B_R(p)$, and let $f:\bb{R} \to \bb{R}$ be a smooth monotone function. Define $\phi:\G(u) \to \bb{R}$ as $\phi:=f \circ u_g$ and assume that $\phi$ is strictly positive. It holds 
    \[
    \m_g(D_{\bar{x},r})^{1/N} \Big(\int_{D_{\bar{x},r}}\phi^{\frac{N}{N-1}} \, d\m_g \Big)^{\frac{N -1}{N}} \leq Cr
    \Big(\int_{D_{\bar{x},r+ \tau}} \lip_g(\phi) \, d\m_g+\frac{1}{\tau}\int_{D_{\bar{x},r}}\phi \, d\m_g \Big)
    \]
    and
    \[ 
    \m_g(D_{\bar{x},r})^{1/N} \Big(\int_{D_{\bar{x},r}}\phi^{\frac{2N}{N-1}} \, d\m_g \Big)^{\frac{N -1}{N}} \leq C
    \Big(r^2\int_{D_{\bar{x},r+ \tau}} \lip_g(\phi)^2 \, d\m_g+
    \frac{2r}{\tau}\int_{D_{\bar{x},r+ \tau}}\phi^2 \, d\m_g \Big).
    \]
\end{thm}

We only sketch the proof of Theorem \ref{T10}, so that one sees how the machinery of the previous section replaces integration by parts in the smooth setting. The part of the proof that we omit is formally the same as the one in \cite{Ding}.

 \begin{thm} \label{T10}
 Let $(\X,\sd,\m)$ be an $\RCD(0,N)$ space, let $B_R(p) \subset \X$ be a ball, and let $u \in \W^{1,1}_{loc}(B_R(p))$ be a positive solution of the minimal surface equation. There exists a constant $C(N)>0$ such that, setting $\bar{p}:=(p,u(p))$, it holds
     \[
         \sup_{\G(u) \cap D_{\bar{p},R/2}} u_g \leq C \inf_{\G(u) \cap D_{\bar{p},R/2}} u_g.
    \]
     \begin{proof}
         Let $0<s<R$ and define $w: \G(u) \to \bb{R}$ by
         \[
         w:=\log{u_g} - \fint_{D_{\bar{p},s}} \log{u_g} \, d\m_g.
         \]
         Let $\eta \in \Lip_c(\G(u))$ be a function to be determined later and let $q \in [0,+ \infty)$. Because of Corollary \ref{CCC10}, it holds
         \begin{align*}
         \int_{\G(u)} & \lip_g(w)^2 \eta^2 |w|^q \, d\m_g=
         \int_{\G(u)} \lip_g(w) \cdot \lip_g(\eta^2 |w|^q) \, d\m_g
         \\
         & =2\int_{\G(u)}\eta |w|^q \lip_g(\eta) \cdot \lip_g(w) \, d\m_g +q\int_{\G(u)} \eta^2|w|^{q-1}w \lip_g(w)^2 \, d\m_g
         \\
         & \leq \frac{1}{2}\int_{\G(u)} \lip_g(w)^2 \eta^2 |w|^q \, d\m_g
         +2\int_{\G(u)} \lip_g(\eta)^2|w|^q \, d\m_g+q\int_{G(u)}\eta^2|w|^{q-1}\lip_g(w)^2 \, d\m_g.
         \end{align*}
         The previous inequalities imply
         \[
         \int_{\G(u)} \lip_g(w)^2 \eta^2 |w|^q \, d\m_g 
         \leq 4 \int_{\G(u)}\lip_g(\eta)^2|w|^q \, d\m_g +2q\int_{\G(u)}\eta^2 |w|^{q-1}\lip_g(w)^2 \, d\m_g,
         \]
         so that choosing $\eta$ appropriately we obtain, for every $r \leq R/4$, 
         \[
         \int_{D_{\bar{p},3r}} \lip_g(w)^2 \, d\m_g \leq \frac{8}{r^2}\m_g(D_{\bar{p},4r}).
         \]
         This corresponds to equation $(4.21)$ in \cite{Ding}. 
         The remaining part of the proof in \cite{Ding} consists of the following steps:
         \begin{itemize}
             \item Let $f \in \sL^\infty_{loc}(\G(u))$, $k \geq 0$ and $r \in (0,R)$. We set
             \[
             \|f\|_{k,r}:=\Big(\fint_{D_{\bar{p},r}} |f|^k \, d\m_g \Big)^{1/k}. 
             \]
             Applying the Neumann-Poincarè inequality \cite[Lemma 3.5]{Ding}, the Sobolev inequality \cite[Lemma 3.4]{Ding} and integration by parts arguments, one deduces \cite[Equation (4.34)]{Ding}: for a sequence of radii $r_j \downarrow r/2$, for some $c(N)>0$, and for every $q \geq 1$, it holds
             \[
             \| |w|^q \|_{\frac{Nq}{N-1}r_{j+2}} \leq (c(N)2^j)^{\frac{1}{q}} (\|w\|_{q,r_j}+4q).
             \]
             In our setting, integration by parts can be carried out as in the begininning of the proof thanks to Theorems \ref{CT5} and \ref{CT6}. The Neumann-Poincarè inequality \cite[Lemma 3.5]{Ding} is replaced by the formally identical Theorem \ref{T8}, while the Sobolev inequality \cite[Lemma 3.4]{Ding} is replaced by Theorem \ref{T9}.
             \item By an iteration argument, this implies, for a different $c(N)>0$, \cite[Equation (4.40)]{Ding}: $\| w \|_{k,r/2} \leq c(N)k$ for every $k \in \bb{N}$. Using Stirling's formula and the Taylor expansion of the exponential, this implies, for some $\lambda(N)>0$, \cite[Equation (4.44)]{Ding}:
             \begin{equation} \label{E65}
             \int_{D_{\bar{p},r/2}} u_g^{\lambda(N)} \, d\m_g \int_{D_{\bar{p},r/2}} u_g^{-\lambda(N)} \, d\m_g \leq c(N).
             \end{equation}
             This part follows from algebraic manipulations and it is formally the same in our setting.
             \item The conclusion follows combining \eqref{E65} and \cite[Theorem 4.2]{Ding}, according to which $u^\lambda$ satisfies a mean value inequality on $\G(u)$.
             \cite[Theorem 4.2]{Ding} is once again obtained combining integration by parts arguments and the Sobolev inequality \cite[Lemma 3.4]{Ding}. The same argument works in our setting,  using Theorem \ref{T9} instead of \cite[Lemma 3.4]{Ding}.
         \end{itemize}
     \end{proof}
 \end{thm}

As an immediate application of the Harnack inequality we get the following result, which implies Theorem \ref{CT1}.
 \begin{thm} \label{T13}
 Let $(\X,\sd,\m)$ be an $\RCD(0,N)$ space and let $u \in \W^{1,1}_{loc}(\X)$ be a function satisfying one of the equivalent conditions of Theorem \ref{T11}. If $u$ is positive, then it is constant.
\end{thm}

\subsection{Applications to the smooth setting}

In this section we prove Theorems \ref{CT2} and \ref{CT3} from the Introduction.
Given a manifold $(\ssf{M},g)$ we denote by $\ssf{Vol}_g$ its volume measure and by $\ssf{d}_g$ its distance. If $V:\ssf{M} \to \bb{R}$ is a smooth function we say that the metric measure space $(\ssf{M}^n,\ssf{d}_g,e^{-V} \ssf{Vol}_g)$ is a weighted manifold.
Given an open set $\Omega \subset \M$ we say that a function $u \in C^{\infty}(\Omega)$ is a solution of the weighted minimal surface equation on $\Omega \setminus \partial \M$ if
\[
\Div \Big( \frac{e^{-V} \nabla u}{\sqrt{1+|\nabla u|^2}} \Big)=0 \quad \text{on } \Omega \setminus \partial \M.
\]
We say that the boundary of a manifold with boundary is convex if its second fundamental form w.r.t. the inward pointing unit normal is positive.  
The next proposition can be obtained repeating an argument in \cite[Theorem $2.4$]{H20}.

\begin{proposition} \label{P29}
Let $(\M^n,\sd_g,e^{-V}  \ssf{Vol}_g)$ be a  weighted manifold with convex boundary such that, for a number $N \geq n$, it holds
    \begin{equation} \label{E43}
    \Ric_\M+\Hess_V-\frac{\nabla V \otimes \nabla V}{N-n} \geq 0 \quad \text{ on } \M \setminus \partial \M,
    \end{equation}
    with the convention that if $N=n$ only constant weights are allowed.
    Then, $(\M^n,\sd_g,e^{-V}  \ssf{Vol}_g)$ is an $\RCD(0,N)$ space.
\end{proposition}

Given a weighted manifold with boundary $(\M^n,\sd_g,e^{-V} \ssf{Vol}_g)$ and a smooth vector field $A \in \ssf{TM}$,
we define the pointwise divergence in the weighted manifold by $\Div_V A:=\Div A-\nabla V \cdot \nabla A$.

\begin{proposition} \label{P31}
    Let $(\M^n,\sd_g,e^{-V} \ssf{Vol}_g)$ be a  weighted manifold with convex boundary satisfying condition \eqref{E43} for $N>0$, and let $u \in C^{\infty}(B_{R}(x))$ be a solution of the weighted minimal surface equation on $B_{R}(x) \setminus \partial \M$ whose gradient vanishes on $\partial \M \cap B_{R}(x)$. Then, $u$ solves the minimal surface equation on $B_{R}(x)$ in distributional sense.
    \begin{proof}
        We first check the statement testing against $C^{\infty}_c(\M)$ functions. Let $\phi \in C^{\infty}_c(\M)$ and let $\nu$ be the outward unit normal on $\partial \M$. 
        Integrating by parts, we obtain
        \begin{equation} \label{E44}
        \int_\M \frac{\nabla u}{W_u} \cdot \nabla \phi \, e^{-V} \, d\ssf{Vol}_g=
        -\int_\M \phi \, \Div_V \Big( \frac{\nabla u}{W_u} \Big) e^{-V} \, d\ssf{Vol}_g +
        \int_{\partial \M}\frac{\phi}{W_u} \nabla u \cdot \nu e^{-V} \, d\aH^{n-1},
        \end{equation}
        where $\aH^{n-1}$ is the Hausdorff measure w.r.t. $\sd_g$.
Since $u$ solves the minimal surface equation on $\M \setminus \partial \M$ we have that $\Div_V \Big( \frac{\nabla u}{W_u}\Big)=0$ pointwise, while the second summand in \eqref{E44} is zero by the condition on $\nabla u$. 
The same is true for $\phi \in \Lip_c(\M)$ via a standard approximation argument.
    \end{proof}
\end{proposition}

\begin{customthm}{2}
Let $(\ssf{M}^n,\ssf{d}_g,e^{-V} \ssf{Vol}_g)$ be a weighted manifold with convex boundary such that there exists $N>n$ satisfying
\[
    \Ric_\M+\Hess_V-\frac{\nabla V \otimes \nabla V}{N-n} \geq 0 \quad \text{ on } \ssf{M} \setminus \partial \ssf{M}.
\]
If $u \in C^{\infty}(\M)$ is a positive solution of the weighted minimal surface equation on $\M \setminus  \partial \M$ whose gradient vanishes on $\partial \M$, then $u$ is constant.
    \begin{proof}
        The weighted manifold $(\M^n,\sd_g,e^{-V} \ssf{Vol}_g)$ is an $\RCD(0,N)$ space by Proposition \ref{P29}, while $u$ solves the minimal surface equation in distributional sense on the weighted space by Proposition \ref{P31}. The conclusion follows by Theorem \ref{CT1}.
    \end{proof}
\end{customthm}

Let $x \in \X$, $r>0$ and $f :B_r(x) \to \bb{R}$. We define 
\begin{equation} \label{E73}
\Osc_{x,r}(f):=\sup \{|f(y)-f(x)|:y \in B_r(x) \}.
\end{equation}

\begin{thm} \label{theorem3stronger}
    Let $N \in (1,+\infty)$. For every $T,t,r>0$ with $T \geq t$ there exists $R(N,T,t,r)>r$ such that if $(\X,\sd,\m,x)$ is a pointed $\RCD(0,N)$ space and $u \in \W^{1,1}_{loc}(B_{R}(x))$ is a function satisfying one of the equivalent conditions of Theorem \ref{T11} such that $\Osc_{x,r}(u) \geq t$, then $\Osc_{x,R}(u) \geq T$.
    \begin{proof}
        Let $(\X,\sd,\m,x)$ be a pointed $\RCD(0,N)$ space, and let $u \in \W^{1,1}_{loc}(B_{R}(x))$ be a solution of the minimal surface equation for some $R>10r$. Set $(x,u(x))=:\bar{x} \in \G(u)$. 
        Theorem \ref{T10}, together with the argument of \cite[Proof of Theorem 3.26]{Beck} (see also \cite[Section 5]{moserosc}), implies that
        \[
        \sup_{\G(u) \cap D_{\bar{x},s}} u_g - \inf_{\G(u) \cap D_{\bar{x},s}} u_g \leq \frac{1}{c} \Big(\frac{s}{S} \Big)^c \Big(\sup_{\G(u) \cap D_{\bar{x},S}} u_g - \inf_{\G(u) \cap D_{\bar{x},S}} u_g \Big),
        \]
        for some $c(N)>0$ and every $0<s<S<R/2$. Hence, setting $S:=R/3$, it holds
        \begin{align*}
        t &\leq \Osc_{x,r}(u) \leq \sup_{\G(u) \cap D_{\bar{x},2r}} u_g - \inf_{\G(u) \cap D_{\bar{x},2r}} u_g \\
        & \leq \frac{1}{c} \Big(\frac{6r}{R} \Big)^c \Big(\sup_{\G(u) \cap D_{\bar{x},R/3}} u_g - \inf_{\G(u) \cap D_{\bar{x},R/3}} u_g \Big) \leq 
        \frac{1}{c} \Big(\frac{6t}{R} \Big)^c \Osc_{x,R}(u).
        \end{align*}
        The statement follows.
    \end{proof}
\end{thm}

Combining with Propositions \ref{P29} and \ref{P31}, we obtain the following corollary.

\begin{corollary} \label{C12}
    Fix $N \in \bb{N}$. For every $T,t,r>0$ with $T \geq t$ there exists $R(N,T,t,r)>r$ such that if $(\M^n,\sd_g,e^{-V} d \ssf{Vol}_g)$ is a  weighted manifold with convex boundary satisfying condition \eqref{E43} for $n<N$ and $u \in C^{\infty}(B_{R}(x))$ is a solution of the weighted minimal surface equation on $B_{R}(x) \setminus \partial \M$ whose gradient vanishes on $\partial \M \cap B_{R}(x)$ such that $\Osc_{x,r}(u) \geq t$, then $\Osc_{x,R}(u) \geq T$.
    \end{corollary}
 
 \footnotesize
\printbibliography
\end{document}